\documentclass[12pt, leqno]{article}
\usepackage{amssymb,amsmath,amsthm}
\usepackage{latexsym}
\usepackage[dvips]{epsfig,graphics}

\newcommand{\Atwo}{\mbox{$\mathbb{Z}[x]/(x^2)$}}

\newcommand\A{{\mathcal A}}
\newcommand\Z{{\mathbb{Z}}}
\newcommand\ga{{\alpha}}

\newcommand\epsi{{\varepsilon}}
\newtheorem{theorem}{Theorem}

\newtheorem{conjecture}[theorem]{Conjecture}
\newtheorem{corollary}[theorem]{Corollary}

\newtheorem{definition}[theorem]{Definition}

\newtheorem{lemma}[theorem]{Lemma}

\newtheorem{proposition}[theorem]{Proposition}
\newtheorem{remark}[theorem]{Remark}

\newtheorem{fact}[theorem]{Fact}
\newtheorem{assumptions on algebras}[theorem]{Assumptions on algebras}
\input{epsf}
\def\stackunder#1#2{\mathrel{\mathop{#2}\limits_{#1}}}

\begin{document}

\title{Torsion in Graph Homology}
\author{Laure Helme-Guizon,  J\'{o}zef H.
Przytycki and Yongwu Rong\footnote{Partially supported by NSF grant
DMS-0513918} }
\maketitle

\begin{abstract}
\noindent Khovanov homology for knots has generated a flurry of
activity in the topology community. This paper studies the Khovanov
type cohomology for graphs with a special attention to torsion. When
the underlying algebra is $\mathbb{Z}[x]/(x^2)$, we determine
precisely those graphs whose cohomology contains torsion. For a
large class of algebras, we show that torsion often occurs. Our
investigation of torsion led to other related general results. Our
computation could potentially be used to predict the
Khovanov-Rozansky $sl(m)$ homology of knots (in particular $(2,n)$
torus knot). We also predict that our work is connected with
Hochschild and Connes cyclic homology of algebras.
\end{abstract}

\tableofcontents

\bigskip
\section{Introduction}
In \cite{K00}, M. Khovanov introduced a graded cohomology theory for
knots whose Euler characteristic is the Jones polynomial of the
knot. Since then, there has been a great deal of attention in his
theory. In \cite{HR04}, a cohomology theory for graphs was
constructed using the algebra $\mathbb{Z}[x]/(x^2)$. In \cite{HR05},
this was generalized by allowing more general algebras. In both
papers, the graded Euler characteristic of the homology groups
yields the chromatic polynomial, through an appropriate change of
variable.

The Euler characteristic is defined in terms of the free part of
the cohomology groups. On the other hand, the torsion of these
homology groups often demonstrate more interesting properties. In
\cite{S04}, A. Shumakovitch did an extensive computation and
proved, among other things, that nontrivial alternating knots
always have $\Z_2$-torsion. In \cite{AP04}, M. Asaeda and J.
Przytycki generalized this to adequate links, and located specific
torsions for a large family of links.

Here, we analyze torsions in the graph cohomology introduced in
\cite{HR04} and \cite{HR05}. For the algebra $\mathbb{Z}[x]/(x^2)$,
we determine precisely those graphs whose cohomology contains
torsion in Theorem \ref{torsion theorem A2}. We also study to what
extent this can be extended for more general algebras.

In Section 2, we review the constructions in \cite{HR04} and
\cite{HR05}. In Section 3, we prove some thickness type theorem for
the cohomology groups. Section 4 contains the Torsion Theorem for
the algebra $\mathbb{Z}[x]/(x^2)$. In Section 5, we consider the
cohomology groups for more general algebras. Special attention is
paid to the algebra $\A_m=\mathbb{Z}[x]/(x^m)$ and the graph $P_v$,
the polygon with $v$ vertices. Subsection 5.3 contains a calculation
of cohomology for $\A=\Z[x]/\mathbb I$, where $\mathbb I$ is an
ideal, in the case of a triangle. We end Section 5 with some general
comments.

\section{The graph cohomology, a quick review}
\label{section graph cohomology} We briefly review our
constructions in \cite{HR04}\cite{HR05}. Recall that a
\emph{graded $\Z$-algebra} $\A$ is a $\Z$-algebra with direct sum
decomposition $\A=\oplus_{j=0}^{\infty} A_j$ into $\Z$-submodules
such that $a_i a_j \in A_{i+j}$ for all $a_i\in A_i$ and $a_j\in
A_j$. The elements of $A_j$ are called \emph{homogeneous} elements
of degree $j$.

Since $\A$ can be considered as a $\Z$-module, its graded
dimension is defined as follows.
Let $M=\oplus_{j=0}^{\infty}M_{j}$ be a graded $\mathbb{Z}$-module
where $M_{j}$ denotes the set of homogeneous elements of degree j
of $M$. Recall that the \emph{free rank} of a $\Z$-module $M$ is
$\dim_{\mathbb{Q}}(M \otimes_{\Z} \mathbb{Q})$ and is denoted by
$rank(M)$. Assume that $rank (M_j) <\infty$ for each $j$. The
\emph{ graded dimension of }$M$ is the power series
\[
q\dim M:=\sum _{j=0}^{\infty} q^{j}\: rank(M_{j}).
\]

We will use the notation $\{ .\}$ to denote the \textit{degree
shift} in graded $\mathbb{Z}$-modules. For example,
$\mathbb{Z}^3\{2\}$ denotes the rank three free abelian group
whose elements are of degree two.

 From now on, we work with
algebras that satisfy the following conditions:
\begin{assumptions on algebras}
\label{assumptions on algebras}
 $\A=\oplus_{i=0}^{\infty} A_i$ is a commutative, graded $\Z$-algebra
with identity such that each $A_i$ is free of finite rank.
\end{assumptions on algebras}
Note that these assumptions can sometimes be relaxed. For instance,
if there is no identity or if the $A_i$'s are not free, the
construction can still be made but some properties will have to be
modified.


Now, we fix an algebra $\A$ satisfying Assumption \ref{assumptions
on algebras}. Let $G$ be a graph and $E=E(G)$ be the edge set of
$G$. Let $n=|E|$ be the cardinality of $E$. We fix an ordering on
$E$ and denote the edges by $e_1, \cdots, e_n$. Consider the
$n$-dimensional cube $\{ 0, 1\}^E = \{ 0, 1\}^n$. Each vertex
$\alpha$ of this cube corresponds to a subset $s=s_{\alpha}$ of
$E$, where $e_i \in s_{\alpha}$ if and only if $\alpha_i=1$. The
\emph{height} $|\alpha|$ of $\alpha$, is defined by $|\alpha|=\sum
\alpha_i$, which is also equal to the number of edges in
$s_{\alpha}$.

For each vertex $\alpha$ of the cube, we associate the graded
$\mathbb{Z}$-module $C^{\alpha}(G)$ as follows. Consider $[G:s]$,
the graph with vertex set $V(G)$ and edge set $s$. We assign a
copy of $\A$ to each component of $[G:s]$ and then taking tensor
product over the components. Let $C^{\alpha}(G)$ be the resulting
graded $\mathbb{Z}$-module, with the induced grading from $\A$.
Therefore, $C^{\alpha}(G) \cong \A^{\otimes k}$ where $k$ is the
number of components of $[G:s]$. We define the $i^{\mbox{\tiny
th}}$ chain group of our complex to be
$C^i(G):=\oplus_{|\alpha|=i} C^{\alpha}(G)$.

Next, we describe the differential maps $d^i: C^i(G)\rightarrow
C^{i+1}(G)$. This will be expressed in terms of the \emph{per edge
maps}. Each edge $\xi$ of $\{ 0, 1\}^E$ can be labeled by a sequence
in $\{ 0, 1, *\}^E$ with exactly one $*$. The tail of the edge,
denoted by $\alpha_1$, is obtained by setting $*=0$; and the head,
denoted by $\alpha_2$, is obtained by setting $*=1$. The
\emph{height} $|\xi|$ is defined to be the height of its tail, which
is also equal to the number of 1's in $\xi$.  For $j=1$ and 2, the
$\Z$-module $C^{\alpha_j}(G)$ is $\A^{\otimes k_j}$ where $k_j$ is
the number of connected components of $[G:s_j]$ (here $s_j$ stands
for $s_{\alpha_j}$). Let $e$ be the edge with $s_2=s_1\cup \{ e\}$.
The per-edge map $d_{\xi}: C^{\alpha_1}(G)\rightarrow
C^{\alpha_2}(G) $ is defined as follows.

If $e$ joins a component of $[G:s_1]$ to itself, then $k_1=k_2$
and the components of $[G:s_1]$ and the components of $[G:s_2]$
naturally correspond to each other. We let $d_{\xi}$ to be the
identity map.

 If $e$ joins two different components of $[G:s_1]$,
say $E_1$ to $E_2$ where $E_1, E_2, \cdots, E_{k_1}$ are the
components of $[G:s_1]$, then $k_2=k_1 - 1$ and the components of
$[G:s_2]$ are $E_1\cup E_2 \cup \{ e\}, E_3, \cdots, E_{k_1}$. We
define $d_{\xi}$ to be the identity map on the tensor factors
coming from $E_3, \cdots, E_{k_1}$, and $d_{\xi}$ on the remaining
tensor factors to be the multiplication map $\A \otimes \A
\rightarrow \A$ sending $x\otimes y$ to $ xy$.

The differential $d^i:C^i(G) \rightarrow C^{i+1}(G)$ is then
defined by $d^i = \sum_{|\xi|=i} (-1)^{\xi} d_{\xi}$, where
$(-1)^{\xi}=-1$ (resp. 1) if the number of 1's in $\xi$ before $*$
is odd (resp. even).

\bigskip
It was proved in \cite{HR05} that this defines a graded cochain
complex $\mathcal{C}$ whose graded cohomology groups $H^i$ are
invariants of the graph (i.e. independent of the ordering of the
edges), and whose graded Euler characteristic is the chromatic
polynomial of $G$ evaluated at $\lambda=q\dim \A$, i.e.
\begin{equation}
\label{chromatic as graded Euler Char}
 \chi_{q}(\mathcal{C})=\sum
_{0\leq i\leq n}(-1)^{i} q\dim (H^{i})=\sum_{0\leq i\leq n} (-1)^{i}
q\dim (C^{i})=P_G(q\dim \A)
\end{equation}

The cochain complex $\mathcal{C}$ splits into a sequence of
cochain complexes, one for each degree $j$. More precisely, let
$\mathcal{C}$ = $0\rightarrow C^{0}\rightarrow C^{1}\rightarrow
...\rightarrow C^{n}\rightarrow 0$ be a graded cochain complex
with a degree preserving differential. Decomposing elements of
each cochain group by degree yields $C^{i}=\oplus_{j\geqslant
0}C^{i,j}$, where $C^{i,j}$ is the set of homogeneous elements of
$C^i$ of degree $j$. Since the differential is degree preserving,
the restriction to elements of degree $j$, i.e. $0\rightarrow
C^{0,j}\rightarrow C^{1,j}\rightarrow ...\rightarrow
C^{n,j}\rightarrow 0$ is a cochain complex denoted by
$\mathcal{C}^{j}$. It is clear that $\mathcal{C}$ is the direct
sum of these cochain complexes.

It is also convenient to use the notion of ``enhanced states".
When the algebra is $\mathbb{Z}[x]/(x^2)$, an \textit{enhanced
state} of a graph $G$ is $S=(s,c)$, where $s\subseteq E(G)$ and
$c$ is an assignment of $1$ or $x$ to each connected component of
the spanning subgraph $[G:s]$.  Such enhanced states form a basis
for the chain groups. Furthermore, the differential map $d$ can be
described explicitly in terms of enhanced states. Details can be
found in \cite{HR04}. This can be generalized to more general
algebras with a given generating set for its additive group, with
proper translations for any possible relations on the generators.

\bigskip
An interesting property of these cohomology groups is a long exact
sequence described below, which corresponds to the
deletion-contraction rule of the chromatic polynomial. Let $e$ be
an edge of the graph $G$. Recall that $G-e$ denotes the graph
obtained by deleting $e$ from $G$, and $G/e$ denotes the graph
obtained by contracting $e$ to a point. We have

\begin{theorem}$\;$\\
\label{long exact sequence}(a) For each $i$, $j$, there is a short
exact sequence of graded chain homomorphisms: $0\rightarrow
C^{i-1,j}(G/e)\overset{\alpha}{\rightarrow} C^{i,j}(G)
\overset{\beta}{\rightarrow} C^{i,j}(G-e) \rightarrow 0$, and
therefore, \\
(b)  it induces a long exact sequence of cohomology groups:
$0\rightarrow H^{0,j}(G)\rightarrow H^{0,j}(G-e) \rightarrow
H^{0,j}(G/e)\rightarrow H^{1,j}(G)\rightarrow \ldots \rightarrow
H^{i,j}(G)\overset{\beta ^{\ast }}{\rightarrow }H^{i,j}(G-e)
\rightarrow H^{i,j}(G/e)\rightarrow\ldots \ $

\end{theorem}
This long exact sequence has proved to be a useful computational
tool, as we will see in the next several sections.

\bigskip

Another useful result from \cite{HR05} is about pendant edges.
Recall that a {\em pendant vertex} in a graph is a vertex of degree
one, and a {\em pendant edge} is an edge connecting a pendant vertex
to another vertex.

\begin{proposition}
\label{pendant edge} Let $\A$ be an algebra satisfying Assumption
\ref{assumptions on algebras}. It can be written $\A =
\mathbb{Z}1_{\A\ } \oplus \A'$ as a $\mathbb{Z}$-module, where
$\mathbb{Z}1_{\A }$ is generated by the identity of $\A$ and $\A'$
is a submodule of $\A$. If $e$ is a pendant edge of $G$, then
$H^i(G)\cong H^i(G/e)\otimes \A'$.
\end{proposition}

\section{General results}
\label{section general results} Theorems \ref{vanishing theorem},
\ref{ThicknessPointed}, \ref{Thickness max degree is m} and
Corollary \ref{Thickness-Am} below indicate which cohomology
groups of a graph may be non-trivial. They generalize similar
results that can be found in \cite{H05} in the case $\A=\A_2$.

\subsection{The Vanishing Theorem}

\begin{theorem}
\label{vanishing theorem} Let $G=G_1 \sqcup G_2$, where $G_2$ is
the union of the components of $G$ that are trees, and $G_1$ is
the union of the remaining components. Let $v_i$ be the number of
vertices of $G_i$ and $\mu_i$ be the number of components of $G_i$
($i=1, 2$).
 For any algebra $\A$
satisfying Assumption \ref{assumptions on algebras}, we have
\begin{enumerate}
\item [(a)] $H_{\A}^i(G)\neq 0 \ \ \Rightarrow  \ \ 0 \leq i \leq v_1 - 2\mu_1$,\\
\item [(b)] $Tor(H_{\A}^i(G))\neq 0 \ \ \Rightarrow  \ \ 1 \leq i
\leq v_1-2\mu_1$.
\end{enumerate}
\end{theorem}
In particular, if a graph has no isolated vertices, we get the
following more convenient inequality.
\begin{corollary}
Let $G$ be a $\mu$-component graph with $v$ vertices. If $G$ has
no isolated vertices, then $H^i_{\A}(G)\neq 0 \ \ \Rightarrow \ \
0 \leq i \leq v - 2\mu$.
\end{corollary}
\begin{proof}
We use the notations of Theorem \ref{vanishing theorem}. We have
$v - 2\mu=v_1 - 2\mu_1 +v_2 - 2\mu_2$. By assumption, each
component of $G_2$ has at least 2 vertices, therefore $v_2 -
2\mu_2 \geq 0$ and this implies $v - 2\mu \geq v_1 - 2\mu_1 $.
\end{proof}

\begin{remark} \label{Remark6} Theorem \ref{vanishing theorem} implies the following
result of \cite{H05}:
 If $G$ is a $\mu$-component graph with $v$ vertices, then $H^i_{\A}(G)\neq 0$ implies $0 \leq i \leq v -
 \mu-1$.\\
 In particular, if $v \geq 2$, $H^i_{\A}(G)\neq 0$ implies $0 \leq i \leq v -
2$.
  This observation was
first made by Michael and Sergei Chmutov \cite{C05} in the case
when the graph is connected and $\A=\A_2$.
\end{remark}

\begin{remark} The inequality in Theorem \ref{vanishing theorem} is sharp for graphs that are forests,
polygons, or their disjoint unions. But if a (non-forest) graph
contains a pendant edge, we can improve the previous result. For
example, we have
\end{remark}

\begin{corollary}[of Theorem \ref{vanishing theorem}]
\label{pendant edge vanishing} \ \\(i) $H^{v_G-k-1}(G)=0$ whenever
$G$ has $k$ pendant edge for
$v_G-k-1 \geq 1$.\\
(ii) If $G$ has a pendant edge in a non-tree component, then
$H^i_{\A}(G)\neq 0 \ \ \Rightarrow  \ \ 0 \leq i \leq v_1 -
2\mu_1-1$.
\end{corollary}

\begin{proof}
(i) Remark 6 implies that for any graph $G$ with $v_G$ vertices,
$v_G\geq 2$, $H^i_{\A}(G)\neq 0$ implies $0 \leq i \leq v_G - 2$.
Using Proposition \ref{pendant edge}, we see that if a graph $Y^p$
is obtained from $Y$ by adding a pendant edge then
$H^{i,*}_{\A}(Y^p)= 0$ whenever $H^{i,*}_{\A}(Y)=0$. Combining these
two results, we get that $H^{v-2,*}_{\A}(Y^p)= 0$, where $v$ is the
number of vertices of $Y^p$. We get the announced result by
induction on the number of pendant edges.

(ii)Let $e$ be a given pendant edge in a non-tree component.  By
Proposition \ref{pendant edge}, we have $H_{\A}^i(G) \cong
H_{\A}^i(G/e) \otimes \A'$. Therefore, $H^i_{\A}(G)\neq 0$ implies
$H_{\A}^i(G/e)\neq 0$ which implies $i\leq v_1(G/e)-2\mu_1(G/e) =
v_1(G) - 2\mu_1(G) -1$.
\end{proof}

\begin{proof} [Proof of Theorem \ref{vanishing theorem}]  (a)  We will prove
the inequality by  induction on $n$, where $n$ is the number of
edges in $G$.

If $n=0$, $G$ is the null graph with $v_2$ vertices.  We have
$v_1=\mu_1=0$. It is well-known that $H^0_{\A}(G)\cong \A^{\otimes
v_2}$ and $H^i_{\A}(G)=0$ for all $i>0$ \cite{HR05}.  Thus the
inequality holds.

Suppose the inequality holds for all graphs with no more than $n$
edges. Now let $G$ be a graph with $n$ edges.

If $G$ is a forest, then it is well-known that $H^i(G)=0$ for all
$i>0$ \cite{HR05}. Thus the inequality $0\leq i \leq v_1-2\mu_1$
becomes $0\leq 0 \leq 0$ which is true.

If $G$ has multiple edges, deleting one such edge $e$ yields a
graph $G-e$ with one less edge. Note that $v_1(G-e)=v_1(G)$ and
$\mu_1(G-e)=\mu_1(G)$. Since cohomology groups remain the same
when deleting a multiple edge, applying induction $G-e$ proves
that the inequality holds for $G$.

If $G$ has a loop, the inequality holds trivially since
$H^i_{\A}(G)=0$ for all $i$.

Now, we assume that $G$ is not a forest, has no loop or multiple
edges. It follows that there is a cycle of length 3 or larger. Let
$e$ be an edge in this cycle. Assume that $i>v_1 - 2\mu_1$.
Applying the exact sequence on $(G, e)$, we have

\begin{equation*}
\cdots \rightarrow H_{\A}^{i-1}(G/e)\rightarrow
H_{\A}^{i}(G)\rightarrow H_{\A}^{i}(G-e)\rightarrow \cdots
\end{equation*}

For the graph $G/e$, we have $v_1(G/e)=v_1(G)-1,
\mu_1(G/e)=\mu_1(G)$, thus $v_1(G/e) - 2\mu_1 (G/e) =
v_1(G)-2\mu_1(G) -  1 < i-1$, which implies $H_{\A}^{i-1}(G/e)=0$
by induction.

For the graph $G-e$, we have $v_1(G-e)=v_1(G)$ and $
\mu_1(G-e)=\mu_1(G)$. Thus $v_1(G-e) - 2\mu_1(G-e)=v_1(G) -
2\mu_1(G)<i$, which implies $H_{\A}^{i}(G-e)=0$.

Therefore, we have $H_{\A}^i(G)=0$ by the exact sequence.

(b). By (a), the only thing we need to show is that $H^0(G)$
contains no torsion. But this follows since $H^0(G) =\ker d^0$,
which is a subgroup of the free abelian group $C^0$.
\end{proof}

\bigskip
 In the next subsection, we show that if we restrict the class of
algebras we deal with, we can refine the result of Theorem
\ref{vanishing theorem}.

\subsection{Thickness of the cohomology}
Indeed, this subsection contains the two following results: Theorem
\ref{ThicknessPointed} shows that under the mild assumption that the
algebra is pointed (see Definition  9 below), the cohomology groups
corresponding to ``small" values $i$ and $j$ are trivial.

Theorem \ref{Thickness max degree is m} shows that under the
assumption that the algebra is of degree less than or equal to
$m$, the cohomology groups corresponding to ``big" values $i$ and
$j$ are trivial.

We then combine these results to determine which cohomology groups
may be non trivial when $\A=\A_m$.

\begin{definition}
We say that a graded $\Z$-algebra is \emph{pointed} if it has a
direct sum decomposition into $\Z$-submodules of the form $\A=\Z1
\oplus \A'$ where all the homogeneous elements of $\A'$ have
degree greater than or equal to $1$.
\end{definition}

\begin{theorem}
Let $G$ be a $\mu$-component graph with $v$ vertices. Let $\A$ be
a pointed $\Z$-algebra satisfying Assumptions \ref{assumptions on
algebras}. We have
\begin{enumerate}
\item [(a)] $H_{\A}^{i,j}(G)\neq 0 \Rightarrow
i+j \geq v-\mu $\\
\item [(b)] $Tor(H_{\A}^{i,j}(G))\neq 0 \Rightarrow
i+j \geq v-\mu+1$
\end{enumerate}\label{ThicknessPointed}
\end{theorem}

\begin{proof}
Let $n$ be the number of edges of $G$. We are going to prove the
result by induction on $n$, where our induction hypothesis consists
of the statements (a) and (b).

\bigskip
\emph{Initial condition step.} If $n=0$, $G$ is $N_{v}$, the null
graph on $v$ vertices. The fact that the theorem holds for any
forest is shown below (without using induction) so we won't repeat
the proof here.

\bigskip
\emph{Induction step.} Assume the induction hypothesis holds for all
$n'< n$.

If $G$ has a loop, the induction hypothesis is satisfied, so we
can assume from now on that $G$ has no loop.

If $G$ has multiple edges, the graph $G'$ obtained from $G$ by
replacing multiple edges by single edges has fewer edges so we can
use the induction hypothesis for $G'$. Since $G$ and $G'$ have the
same cohomology groups, the same number of vertices and the same
number of components, this proves the result for $G$. From now on
in this proof, we assume that $G$ has no loop and no multiple
edges.

\bigskip
 \emph{Case 1: $G$ is a forest.}
Assume $G$ is a forest made of $\mu$ trees for a total of $n$ edges
and denote it $F_{\mu,n}$. Note that $F_{\mu,n}$ has $\mu +n$
vertices. The cohomology groups of $F_{\mu,n}$ are obtained from the
cohomology groups of $N_{\mu}$, the graph with $\mu$ vertex and no
edges, by adding $n$ pendant edges. By Theorem \ref{pendant edge},
we get that the cohomology groups are trivial except for
$H_{\A}^0(F_{\mu,n}) \cong \A^{\otimes \mu} \otimes (\A' )^{\otimes
n}$ where $\A'$ satisfies $\A = \Z 1 \oplus \A'$ with the degree of
all the homogeneous elements of $\A'$ greater than or equal to $1$
(since $\A$ is pointed). Hence the degrees $j$ of the homogeneous
elements of $H_{\A}^0(F_{\mu,n})$ are greater than or equal to $n$.
Substituting $n=v-\mu$ in the above yields $0+j \geq v-\mu $.

Also, all the other cohomology groups are trivial and there is no
torsion so the induction hypothesis is satisfied (without using
induction).

\bigskip
\emph{Case 2: $G$  is a not a forest.} Since we assumed that $G$ has
no loop and no multiple edges, $G$ is a not forest means it contains
a cycle of order $\geq 3.$ Let $e$ be an edge on this cycle. The
edge $e$ is not an isthmus so $G-e$ also has $\mu$ components. Of
course, $G/e$ has $\mu$ components no matter whether $e$ is not an
isthmus or not. Therefore, $\mu=\mu_{G}=\mu_{G-e}=\mu_{G/e}$.

The exact sequence on $(G,e)$ yields
\begin{equation*}
\cdots \rightarrow H_{\A}^{i-1,j}(G-e)\rightarrow
H_{\A}^{i-1,j}(G/e)\rightarrow H_{\A}^{i,j}(G)\rightarrow
H_{\A}^{i,j}(G-e)\rightarrow \cdots
\end{equation*}

We don't have to worry about having $i-1\geq 0$ if we define the
cohomology groups with negative heights to be the trivial group.

Since $G-e$ and $G/e$ have one less edge than $G$, the induction
hypothesis applies to them.
\bigskip

 $\blacktriangle$ Assume that $i+j < v_G-\mu_G$. We need to show that
 $H_{\A}^{i,j}(G)=0$.

The exact sequence on $(G,e)$ gives

\begin{equation*}
\cdots \rightarrow  H_{\A}^{i-1,j}(G/e)\rightarrow
H_{\A}^{i,j}(G)\rightarrow H_{\A}^{i,j}(G-e)\rightarrow \cdots
\end{equation*}
The inequality $i+j < v_{G}-\mu_{G}=v_{G-e}-\mu_{G-e}$ implies
$H_{\A}^{i,j}(G-e)=0$ by part (a) of the induction hypothesis. Also
$i+j < v_{G}-\mu_{G} $ implies $(i-1)+j <
(v_{G}-1)-\mu_{G}=v_{G/e}-\mu_{G/e}$ so $H_{\A}^{i-1,j}(G/e)=0$ by
part (a) of the induction hypothesis. Therefore the sequence
$0\rightarrow H_{\A}^{i,j}(G)\rightarrow 0$ is exact so
$H_{\A}^{i,j}(G)=0$.

\bigskip

$\blacktriangle$ Assume that $i+j < v_{G}+1-\mu_{G} $. We need to
show that $Tor(H_{\A}^{i,j}(G))=0$.

The exact sequence on $(G,e)$ gives

\begin{equation*}
\cdots \rightarrow H_{\A}^{i-1,j}(G-e)\rightarrow
H_{\A}^{i-1,j}(G/e)\rightarrow H_{\A}^{i,j}(G)\rightarrow
H_{\A}^{i,j}(G-e)\rightarrow \cdots
\end{equation*}
The inequality $i+j<v_{G}+1-\mu_{G}$ implies
$(i-1)+j<v_{G}-\mu_{G}=v_{G-e}-\mu_{G-e}$ so $H_{\A}^{i-1,j}(G-e)=0$
by part (a) of the induction hypothesis. The inequality
$i+j<v_{G}+1-\mu_{G}=v_{G-e}+1-\mu_{G-e}$ implies
$H_{\A}^{i,j}(G-e)$ is torsion free by part (b) of the induction
hypothesis. Also $i+j < v_{G}+1-\mu_{G}$ implies $(i-1)+j <
v_{G}-\mu_{G}=v_{G/e}+1-\mu_{G/e}$ so $H_{\A}^{i-1,j}(G/e)$ is
torsion free by part (b) of the induction hypothesis. Therefore the
exact sequence is of the form
\begin{equation*} \cdots \rightarrow
0\rightarrow F\rightarrow H_{\A}^{i,j}(G)\rightarrow F'\rightarrow
\cdots
\end{equation*}
where F and F' are torsion free abelian groups. Note that we are
working with finitely generated abelian groups so free and torsion
free are equivalent notions.

The following claim shows that this is enough to prove that
$H_{\A}^{i,j}(G)$ is free.

\textbf{Claim:} If $H$ is an abelian group that satisfies an exact
sequence of the form
\begin{equation*} 0\rightarrow F\rightarrow
H \overset{h}{\rightarrow} F'\rightarrow \cdots
\end{equation*}where $F$ and $F'$ are free abelian groups, then $H$ is
free.

\textbf{Proof of the claim:} The above exact sequence induces the
following short exact sequence
\begin{equation*} 0\rightarrow F\rightarrow
H \overset{h}{\rightarrow} Im \; h\rightarrow 0
\end{equation*}
Such a sequence splits because $Im \; h$ is free. Therefore, $H\cong
F \oplus Im \; h$, which proves that $H$ is free.
\end{proof}

\bigskip

 We now replace the hypothesis that the algebra is pointed by the hypothesis
  that the algebra is of degree less than or equal to $m$.

\begin{definition}
We say that a graded $\Z$-algebra is \emph{of degree less than or
equal to $m$} if  all the homogeneous elements of $\A$ have degree
less than or equal to $m$.
\end{definition}

\begin{theorem}
\label{Thickness max degree is m} Let $G$ be a graph, not
necessarily connected, with $v$ vertices. Let $\A$ be a $\Z$-algebra
of degree less than or equal to $m-1$ satisfying Assumptions
\ref{assumptions on algebras}. We have $H_{\A}^{i,j}(G)\neq 0
\Rightarrow (m-1)i+j \leq (m-1)v$.
\end{theorem}

\begin{proof}
Let $n$ be the number of edges of $G$. We are going to prove the
result by induction on $n$.

\bigskip
\emph{Initial condition step:} If $n=0$, $G$ is $N_{v}$, the null
graph on $v$ vertices. In particular, $G$ is a forest. The fact that
the theorem holds for any forest is shown below (without using
induction) so we won't repeat the proof here.

\bigskip
\emph{Induction step:} Let $n \geq 1$. Assume the induction
hypothesis holds for $n-1$.

 Assume that $(m-1)i+j>(m-1)v_G$. We need to show that $H_{\A}^{i,j}(G)=0$.

Let $e$ be an edge of $G$. The exact sequence on $(G,e)$ gives
\begin{equation*}
\cdots \rightarrow  H_{\A}^{i-1,j}(G/e)\rightarrow
H_{\A}^{i,j}(G)\rightarrow H_{\A}^{i,j}(G-e)\rightarrow \cdots
\end{equation*}
Note that since $G-e$ and $G/e$ have one less edge than $G$, the
induction hypothesis applies to them. The inequality
$(m-1)i+j>(m-1)v_G=(m-1)v_{G-e}$ implies $H_{\A}^{i,j}(G-e)=0$ by
 the induction hypothesis. Also $(m-1)i+j>(m-1)v_{G}$
implies $(m-1)(i-1)+j>(m-1)(v_{G}-1)=(m-1)v_{G/e}$ so
$H^{i-1,j}(G/e)=0$ by the induction hypothesis. The sequence
$0\rightarrow H_{\A}^{i,j}(G)\rightarrow 0$ is exact so
$H_{\A}^{i,j}(G)=0$.
\end{proof}


 We define the $\Z$-algebra $\A_m$ by $\A_m:=\Z[x]/(x^m)$. In the particular case of an algebra of the form
  $\A_m$ we can combine the previous theorems and get
\begin{corollary}
 \label{Thickness-Am}
Let $G$ be a $\mu$-component graph with $v$ vertices. If $G$ has no
isolated vertices.
\begin{enumerate}
\item $H_{\A_m}^{i,j}(G)\neq 0 \Rightarrow \left\{
\begin{array}{ll}
0 \leq i \leq v-2\mu & (1a)\\
i+j \geq v-\mu & (1b)\\
(m-1)i+j \leq (m-1)v & (1c)
\end{array}
\right.$\\
\item $Tor(H_{\A_m}^{i,j}(G))\neq 0 \Rightarrow \left\{
\begin{array}{ll}
1 \leq i \leq v-2\mu & (2a)\\
i+j \geq v+1-\mu & (2b)\\
(m-1)i+j \leq (m-1)v & (2c)
\end{array}
\right.$
\end{enumerate}
\end{corollary}

This shows the thickness of the region where $H^{i,j}(G)$ is
supported. An illustration is shown in Figure  \ref{Thickness
figure}. In the case when m=2, and G is connected, this was also
observed by M. Chmutov and S. Chmutov \cite{C05}.

\begin{figure}[h]
\begin{center}
\scalebox{.8}{\includegraphics{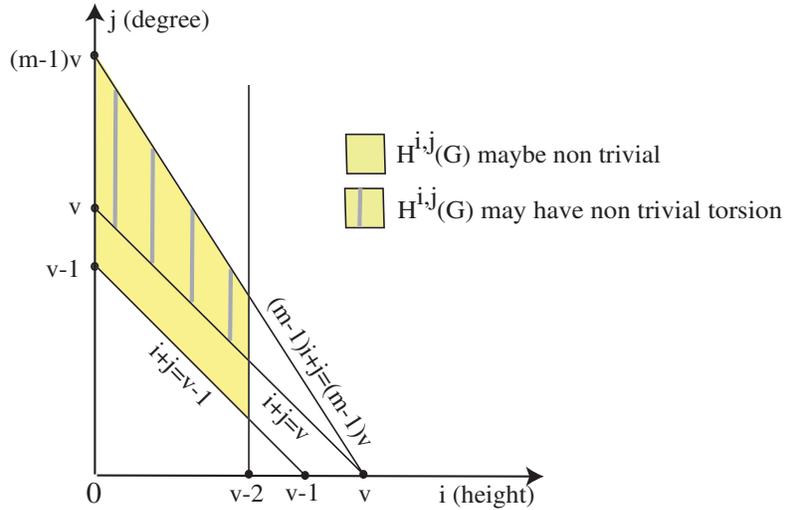}}
\caption{\textit{Thickness restrictions when $\A=\A_m$ for a
connected graph.}} \label{Thickness figure}
\end{center}
\end{figure}

\subsection{The particular case of polygons}

For $n \geq 3$, we denote by $P_n$ the polygonal graph with $n$
edges (and $n$ vertices). We can extend this definition in a natural
way to the ``degenerated" cases $n=1$ and $n=2$ the following way.
Let $P_1$ be the graph made of one vertex with a loop attached to it
and let $P_2$ be the graph made of two vertices connected by two
edges.
\begin{lemma}
Let $P_n$ be the polygon with $n$ edges and let $\A$ be an algebra
satisfying Assumption \ref{assumptions on algebras}.

\[
H_{\A}^{i+1}(P_{n+1}) =H_{\A}^{i}(P_n) \:\text{if } \; i\geq 1 \;
\text{and  } \;n\geq 1.
\]
\label{Polygon induction}
\end{lemma}

\begin{proof}
In \cite{HR05}, it is shown that an algebra $\A$ satisfying the
Assumption \ref{assumptions on algebras} can be written $\A=\Z 1
\oplus \A'$. With these notations, the cohomology groups for a tree
$T_n$ with $n$ edges are known to be (see \cite{HR05}):
\[
\left\{
\begin{array}{l}
H_{\A}^{0}(T_n)=\A \otimes \A'^{\otimes n}\\
H_{\A}^{i}(T_n)=0 \text{ if } i\neq 0.
\end{array}
\right.
\]

Fix $n \geq 1$, fix $i \geq 1$ and let $e$ be an edge of $P_{n+1}$.
Applying the above result for trees to the exact sequence on
$(P_{n+1},e)$ yields:

 $
\begin{array}{ccccccccccc}
&  & G-e &  & G/e &  & G &  & G-e &  &\\
& \cdots \rightarrow  & H^{i}(T_n) & \rightarrow & H^{i}(P_n) &
\rightarrow  & H^{i+1}(P_{n+1}) & \rightarrow & H^{i+1}(T_n)&
\rightarrow  &\cdots  \\
&  & \shortparallel &  &  &  &  &  & \shortparallel  &  &\\
&  & 0 &  &  &  &  &  & 0  &  &\\
\end{array}
$

This proves the expected result.
\end{proof}
By induction, we get
\begin{corollary}
\label{Polygon in terms of H1} Let $P_n$ be the polygon with $n$
edges and let $\A$ be an algebra satisfying that assumptions
(\ref{assumptions on algebras}).
\[
H_{\A}^{i}(P_n)= \left\{
\begin{array}{ll}
H_{\A}^{1}(P_{n-i+1})& \text{if } \; i\leq n,\\
0& \text{otherwise.}
\end{array}
\right.
\]
\end{corollary}
This shows that determining the cohomology groups at height $1$ is
enough to recover all of them. Indeed, the only cohomology groups
not determined by the above formula are the ones at height $0$.
They can be obtained from the previous ones via the chromatic
polynomial since $\sum _{0\leq i\leq n}(-1)^{i} q\dim
(H^{i})=P_G(q\dim \A)$.

\section{Torsion when the algebra is $\Atwo$ }
\label{section algebra is A2} In this section, we will focus on the
algebra $\A_2 = \mathbb{Z}[x]/ (x^2)$. Thus all cohomology groups
will be based on $\A_2$. In 4.1, we review some useful facts. Our
main result is stated in 4.2, and proved in 4.3 and 4.4.

\subsection{Facts and observations}
Corollary \ref{Thickness-Am} specialized to $m=2$ gives us some
information about the thickness of the cohomology of a connected
graph:
\begin{corollary}
\label{Thickness-A2}  The cohomology of a connected graph is
concentrated along two diagonals, namely $i+j=v$ and $i+j=v-1$.
Torsion can occur on one diagonal only, namely $i+j=v$.
\end{corollary}

\begin{figure}[h]
\begin{center}
\scalebox{.7}{\includegraphics{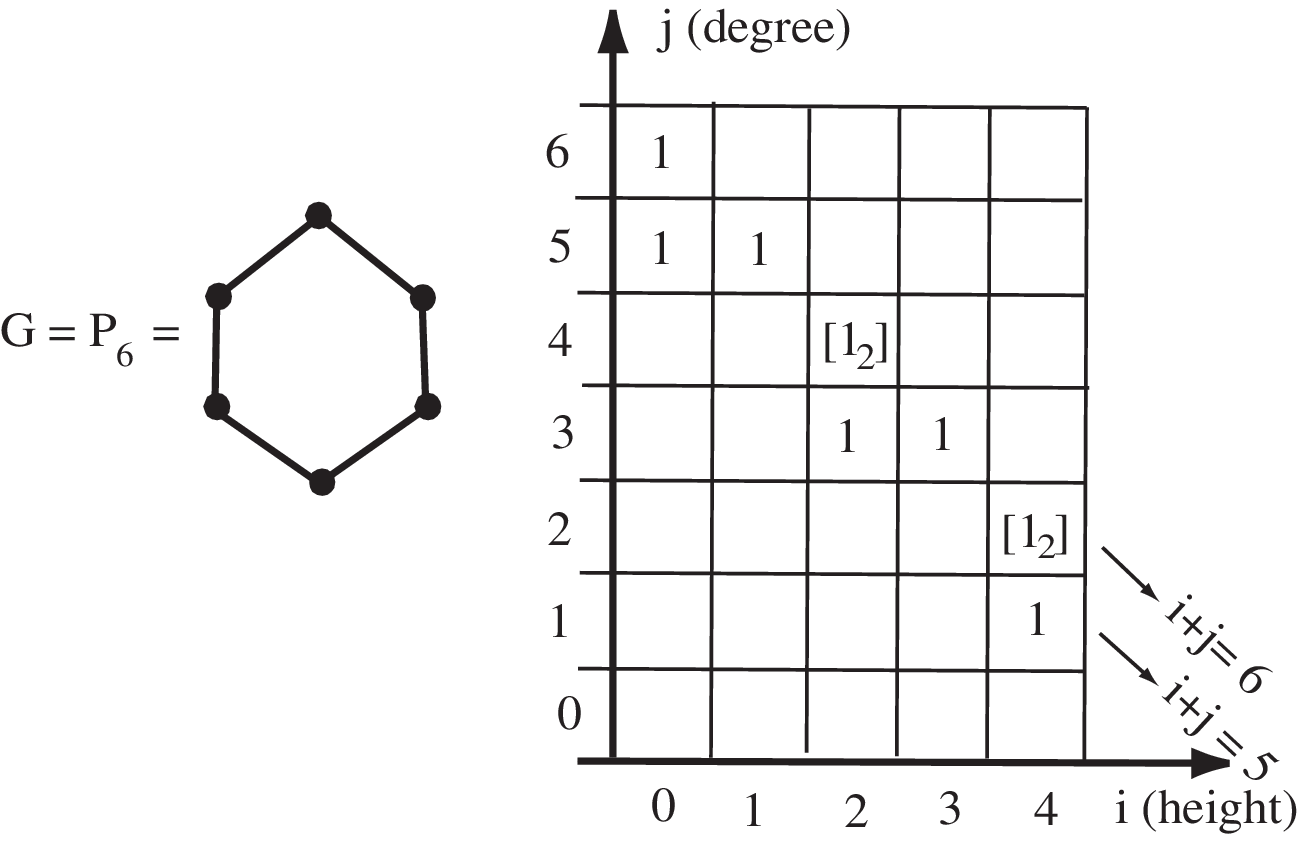}}
\end{center}
\caption{\textit{$P_6$ summary when the algebra is $\Atwo$}}
\label{P6summaryDiag}
\end{figure}

The ``diagonals" language refers to the way of keeping track of
the cohomology groups illustrated in Figure \ref{P6summaryDiag}.
In the array that keeps track of the cohomology groups, the
numbers without brackets indicate the number of copies of
$\mathbb{Z}$ while a number with brackets of the form $[k_{\ell}]$
indicates $k$ copies of $\Z_{\ell}.$ For instance, in the case of
the cyclic graph $P_6$ (see Figure \ref{P6summaryDiag} below), the
$[1_2]$ in position $i=2, j=4$
 means $H^{2,4}(P_6) \cong \mathbb{Z}_2 \{4\}$ and the $1$ in position $i=2, j=$
 means $H^{2,3}(P_6) \cong \mathbb{Z} \{3\}$,
 so $H^{2}(P_6) \cong \mathbb{Z} \{3\}\oplus \mathbb{Z}_2 \{4\}.$

This example also illustrates that all the non trivial cohomology
groups have height $i$ such that $0 \leqslant i \leqslant 6-2$.

\medskip \noindent
We now remind some of the facts that were proved in \cite{HR04} and
that are going to be useful to determine for which graphs some
torsion occur.

\begin{fact}
\label{fact forests} The cohomology groups of a forest, so in
particular the ones of a tree, are all supported on $H^0$ and don't
have torsion.
\end{fact}

\begin{fact}
\label{fact graphs with loops} If the graph has a loop then all the
cohomology groups are trivial.
\end{fact}

\begin{fact}
\label{fact graphs with multiple edges} The cohomology groups are
unchanged if the multiple edges of a graph are replaced by single
edges.
\end{fact}

\medskip \noindent The table below from \cite{HR04} illustrates our computational results
 (up to $n=6$ and $ i=4 $) for circuit graphs,
where $n$ is the length of the cycle and $i$ is the height of the
homology group. We denote by $P_n$ the circuit graph on $n$ edges
where $P$ stands for polygon because the notation $C_n$ has
already been used for the chain groups.

\begin{center}
\begin{tabular}{|c|l|l|l|l|l|}
\hline $n \backslash i$ & $H^0$ & $H^1$ & $H^2$ & $H^3$ & $H^4$ \\
\hline $P_1$ & 0 & 0 & 0 & 0 & 0 \\
\hline $P_2$ &
$\mathbb{Z}\{2\}\oplus \mathbb{Z}\{1\}$ & 0 & 0 & 0 & 0 \\
\hline $P_3$ & $\mathbb{Z}\{3\}$ & $\mathbb{Z}_2\{2\} \oplus
\mathbb{Z}\{1\}$ & 0 & 0 & 0
\\
\hline $P_4$ & $\mathbb{Z}\{4\}\oplus \mathbb{Z}\{3\}$ &
$\mathbb{Z}\{3\} $ & $\mathbb{Z }_2\{2\}\oplus \mathbb{Z}\{1\}$ & 0
& 0 \\
\hline $P_5$ & $\mathbb{Z}\{5\}$ & $\mathbb{Z}_2\{4\}\oplus
\mathbb{Z}\{3\}$ &
$\mathbb{ Z}\{3\}$ & $\mathbb{Z}_2\{2\}\oplus \mathbb{Z}\{1\}$ & 0 \\
\hline $P_6$ & $\mathbb{Z}\{6\}\oplus \mathbb{Z}\{5\}$ &
$\mathbb{Z}\{5\} $ & $\mathbb{Z }_2\{4\}\oplus \mathbb{Z}\{3\}$ &
$\mathbb{Z}\{3\}$ & $\mathbb{Z}_2\{2\}\oplus \mathbb{Z}\{1\}$
\\ \hline
\end{tabular}
\end{center}

Note that this table illustrates Lemma \ref{Polygon induction} and
corollary \ref{Thickness-A2}.

This table illustrates the following general result that appeared in
\cite{HR04} and that we repeat here for convenience.
\begin{theorem}
\label{Homology polygon A2}
\begin{equation*}
\text{For }i>0,H_{\A_2}^{i}(P_{n})\cong \left\{
\begin{array}{ll}
\mathbb{Z}_{2}\{n-i\}\oplus \mathbb{Z}\{n-i-1\} & \text{if
}n-i\geqslant 2
\text{ and is even} \\
\mathbb{Z}\{n-i\} & \text{if }n-i\geqslant 2\text{ and is odd} \\
0 & \text{if }n-i\leqslant 1.
\end{array}
\right.
\end{equation*}
\end{theorem}

\begin{equation*}
\text{For }i=0,H_{\A_2}^{0}(P_{n})\cong \left\{
\begin{array}{ll}
\Z\{n\}\oplus \mathbb{Z}\{n-1\} & \text{if }n\text{ is even and }
n\geqslant 2 \\
\Z\{n\} & \text{if }n\text{ is odd and }n\geqslant 2 \\
0 & \text{if }n=1.
\end{array}
\right.
\end{equation*}

We first note that, for all $n\geq 3$, $H^{\ast }(P_{n})$ contains
torsion.

A closer look at these examples reveals that, in the case of an odd
cycle, there seems to always be a torsion in $H^1(P_v)$ in degree
$v-1$, and in the case of an even cycle, there seems to always be a
torsion in $H^2(P_v)$ in degree $v-2$ where $v$ is the number of
vertices of the graph. These remarks will guide us for the
formulation of the lemmas.

\subsection{The result}

\begin{theorem}
\label{torsion theorem A2} The cohomology $H_{\A_2}^*(G)$ of a graph
$G$ contains a torsion part if and only if  $G$ has no loops and
contains a cycle of order greater than or equal to $3$. In this
case,  if the cycle has odd length $H^{1,v-1}_{\A_2}(G)$ contains a
$\mathbb{Z}_2$-torsion, and if the cycle has even length
$H^{2,v-2}_{\A_2}(G)$ contains a $\mathbb{Z}_2$-torsion.
\end{theorem}

\begin{proof}
``$\Rightarrow$'' If $G$ has a loop, all the cohomology groups are
zero so there can be no torsion. So we can assume $G$ has no
loops. If $G$ does not have a cycle of order greater than or equal
to $3$, then the graph obtained from $G$ by replacing multiple
edges by single edges is a forest and therefore has no torsion in
its cohomology \cite{HR04}. Since deleting multiple edges doesn't
change the cohomology groups \cite{HR04}, the cohomology of $G$
has no torsion either.

 ``$\Leftarrow$''  If $G$ has  no loop and contains a cycle of order greater than or equal to $3$.
 The following two lemmas  \ref{odd cycle case}
and \ref{even cycle case} imply that $H_{\A_2}^*(G)$ contains
torsion.
\end{proof}

\subsection{The odd cycle case}
\begin{lemma}
\label{odd cycle case}If a loopless graph $G$ with $v$ vertices
contains an odd cycle of length $\geq 3$, then $H^{1,v-1}(G)$
contains a $\mathbb{Z}_2$-torsion.
\end{lemma}

\begin{proof}
Let $G$ be a loopless graph with $v$ vertices containing a cycle of
length $2s+1$ with $s\geq 1$.

It suffices to find an element $z$ in $\ker d^1$ of degree $v-1$
such that  $2z=0$  in $H^1(G)$ and $z \neq 0$ in $H^1(G)$. The
condition $2z=0$ in $H^1(G)$ is the same as $2z \in \mbox{Im}
d^{0,v-1}$. Therefore our first step will be to determine $\mbox{Im}
d^{0,v-1}$.

\noindent $\blacktriangle$ Matrix representation of $d^{0,v-1}$:

\noindent For convenience, we label the vertices of the graph
starting with the ones in the cycle. The vertices in the cycle are
labelled monotonically $v_1$ to $v_{2s+1}$, with the requirement
that each $v_i$ is adjacent to $v_{i+1}$ and $v_{2s+1}$ is adjacent
to $v_{1}$. The vertices that are not in the cycle are labelled
$v_{2s+2}$ to $v_{p}$. Examples are given in Figure \ref{odd cycle
example}.

We label the edges in the cycle so that $e_i$ is the edge $v_{i}
v_{i+1}$ for $1 \leq i \leq 2s$ and $e_{2s+1}$ is the edge
$v_{2s+1} v_1$. The edges that are not in the cycle are labelled
$e_{2s+2}$ to $e_{n}$.

We use the notion of enhanced states here. The basis elements of
$C^{0,v-1}(G)$ are enhanced states $(s, c)$ where $s=\emptyset$ and
 $c$ is an assignment of $1$ or $x$ to the $v$ vertices of $G$.
 Since the total degree is $v-1$, all
vertices are assigned the value $x$ except one which is assigned the
value $1$. The basis element for which $1$ is assigned to the vertex
$v_i$ and with $x$ assigned to all the other vertices is denoted by
$b_i$. The basis elements for which the $1$ is assigned to a vertex
in the cycle are $b_1$ to $b_{2s+1}$.

In order to write a matrix for $d^{0,v-1}$, we also need to describe
the basis elements of the target space $C^{1,v-1}(G)$.

First, note that each of these basis elements contains one edge,
which, by assumption, is not a loop. Therefore, each of these basis
elements has $v-1$ components. Thus, for degree reasons, all the
components are assigned the value $x$. The basis element for which
the present edge is $e_i$ and with $x$ assigned to all the
components is denoted by $a_i$. The basis elements for which the
present edge is in the cycle are $a_1$ to $a_{2s+1}$.


\begin{figure}[h]
\centerline{\psfig{figure=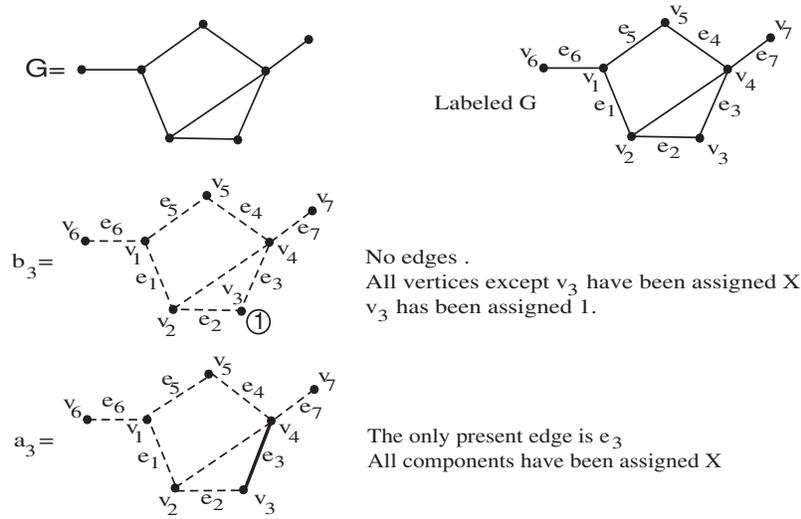,height=7.0cm}}
\caption{\textit{Notation for basis elements basis element in
$C^{0,v-1}(G)$ and $C^{1,v-1}(G)$}} \label{odd cycle example}
\end{figure}

With these notations, we get a matrix representing $d^{0,v-1}$,
which is shown in (\ref{Matrix d0}) below .

\bigskip

\noindent $\blacktriangle$ Let $z=\sum_{i=1}^{n} a_i$

\noindent $\triangle$  We first prove that $2z \in \mbox{Im}
d^{0,v-1}$:

In the matrix representation of $d^{0,v-1}$, the coordinates of the
image of $b_i$ under the differential, in the basis $(a_i)_i$ are
listed in the $i^{\mbox{\tiny th}}$ row.

\begin{equation}
\label{Matrix d0}
\begin{array}{ccc}
\text{\emph{{\small Matrix of }}}d^{0,v-1} &
\begin{array}{cccccccccc}
 & a_{1} & a_{2} & \cdots  & a_{2s+1} & | & a_{2s+2} & \cdots  &  a_{n} &
\end{array}
& \text{Coeff} \\
\begin{array}{c}
d(b_{1})\rightarrow  \\
d(b_{2})\rightarrow  \\
\vdots  \\
\vdots  \\
d(b_{2s+1})\rightarrow  \\
---- \\
d(b_{2s+2})\rightarrow  \\
\vdots  \\
d(b_{v })\rightarrow
\end{array}
& \left[
\begin{array}{cccccccccc}
1 &  &  &  & 1 & | &  &  &  &  \\
1& 1 & &  &  & | &  &  &  &  \\
& \ddots & \ddots  &  &  & | &  & M_{1} &  &  \\
&  & 1 & 1 &  & | &  &  &  &  \\
 &  &  & 1 & 1 & | &  &  &  &  \\
- & - & - & - & - & + & - & - & - & - \\
&  &  &  &  & | &  &  &  &  \\
&  & 0 &  &  & | &  & M_{2} &  &  \\
&  &  &  &  & | &  &  &  &
\end{array}
\right]  &
\begin{array}{c}
1 \\
1 \\
\vdots  \\
1 \\
1 \\
-- \\
1 \\
\vdots  \\
1
\end{array}
\\
\overline{
\begin{array}{c}
d\left( \sum b_{i}\right)
\end{array}
} & \overline{\underbrace{
\begin{array}{cccccccccc}
2, & 2, & \cdots  & \cdots  & 2 & | & 2, & \cdots  & \cdots  & ,2
\end{array}
}} &  \\
& =2z &
\end{array}
\end{equation}

The coefficients (all equal to $1$ in (\ref{Matrix d0})) in the
column at the right indicate the coefficient by which each line is
multiplied before addition.

By adding all the rows of this matrix, we see $d\left(
\sum_{i=1}^{v} b_{i}\right)=2z$ so $2z \in \mbox{Im} d^{0,v-1}$.
This also implies that $z \in \ker d^{1,v-1}$.
 Indeed,
 $2d^1(z)=d^1(2z)=0$ in $C^2(G)$ which implies $d^1(z)=0$ since $C^2(G)$ has no torsion.

The reason why adding the rows of $M_1$ and $M_2$ always yields a
coordinate of $2$ on the $j$th column for $j \geq 2s+2$ is as
follows. Each edge has two ends (remember that there are no loops),
so each of these $\left\{ a_{j}\right\}_{j \geq 2s+2}$ is in the
image of exactly two $b_i$'s, the ones corresponding to a $1$ placed
at each endpoint of the edge $e_{j}$.

\bigskip \noindent $\triangle$ It remains to show that  $z
\notin \mbox{Im} d^{0,v-1}$. Assume $z=d(x)$ for some $x$ in
$C^{0,v-1}$. This $x$ can be written $x=\sum_{i=1}^{v} \alpha_i
b_{i}$ for some $\alpha_i$.

$ \hspace{-1cm}
\begin{array}{ccc}
\text{\emph{{\small Matrix of }}}d^{0,v-1} &
\begin{array}{cccccccccc}
 & a_{1} & a_{2} & \cdots  & a_{2s+1} & | & a_{2s+2} & \cdots  &  a_{n} &
\end{array}
& \text{Coeff} \\
\begin{array}{c}
d(b_{1})\rightarrow  \\
d(b_{2})\rightarrow  \\
\vdots  \\
\vdots  \\
d(b_{2s+1})\rightarrow  \\
---- \\
d(b_{2s+2})\rightarrow  \\
\vdots  \\
d(b_{v})\rightarrow
\end{array}
& \left[
\begin{array}{cccccccccc}
1 &  &  &  & 1 & | &  &  &  &  \\
1& 1 & &  &  & | &  &  &  &  \\
& \ddots & \ddots  &  &  & | &  & M_{1} &  &  \\
&  & 1 & 1 &  & | &  &  &  &  \\
 &  &  & 1 & 1 & | &  &  &  &  \\
- & - & - & - & - & + & - & - & - & - \\
&  &  &  &  & | &  &  &  &  \\
&  & 0 &  &  & | &  & M_{2} &  &  \\
&  &  &  &  & | &  &  &  &
\end{array}
\right]  &
\begin{array}{c}
\alpha _{1} \\
\alpha _{2} \\
\vdots  \\
\alpha _{2s} \\
\alpha _{2s+1} \\
-- \\
\alpha _{2s+2} \\
\vdots  \\
\alpha _{v}
\end{array}
\\
\overline{
\begin{array}{c}
z=d\left( \sum_{i=1}^{v} \ga_i b_{i}\right)
\end{array}
} & \overline{\underbrace{
\begin{array}{cccccccccc}
1, & 1, & \cdots  & \cdots  & 1 & | & 1, & \cdots  & \cdots  & 1
\end{array}
}} &  \\
& =z &
\end{array}
$

This means that the result of multiplying the first line by
$\ga_1$, the second line by $\ga_2$, etc, and adding all the lines
yields $(1,1, \cdots,1)$, the coordinates of $z$ on the $a_i$'s.

\bigskip
If we now read this by columns, we get a contradiction:

$
\begin{array}{ccc}
\begin{array}{c}
\text{ \ Column 1:} \\
\text{+Column 2:} \\
\vdots  \\
\\
\text{+Column 2s+1:}
\end{array}
&
\begin{array}{ccccccccc}
\alpha _{1} & + & \alpha _{2} &  &  &  &  &  &  \\
&  & \alpha _{2} & + & \alpha _{3} &  &  &  &  \\
&  &  &  & \ddots & \ddots &  &  &  \\
&  &  &  &  &  & \alpha _{2s} & + & \alpha _{2s+1}\\
\alpha _{1} & + &  &  &  &  &  & + & \alpha _{2s+1}
\end{array}
&
\begin{array}{c}
=1 \\
=1 \\
\vdots  \\
\vdots  \\
=1
\end{array}
\\
& \overline{\underbrace{
\begin{array}{ccccccccc}
2(\alpha _{1} & + & \alpha _{2} &  &  &  & \alpha _{2s} & + &
\alpha _{2s+1})
\end{array}
}} & \overline{=\underbrace{2s+1}} \\
& \text{even} & \text{odd}
\end{array}
$

This shows that $z \notin \mbox{Im} d^{0,v-1}$.
\bigskip

\end{proof}
\subsection{The even cycle case}

The simplest simple graph without an odd cycle is the ``square"
$C_4$. As mentioned earlier, in this case the torsion appears only
for $H^2(G)$, which indicates that we have to look deeper into the
cohomology to find torsion than in the odd cycle case.

\begin{lemma}
\label{even cycle case} Let $G$ be a simple graph, i.e. a graph with
no loops and no multiple edges. If $G$ contains an even cycle of
length $\geq 4$, then $H^{2,v-2}(G)$ contains a
$\mathbb{Z}_2$-torsion.
\end{lemma}

\begin{proof}
It suffices to find an element $z$ in $\ker d^{2,v-2}$ such that
$2z=0$ in $H^2(G)$ and $z \neq 0$ in $H^2(G)$. The condition
\noindent $2z=0$ in $H^2(G)$ means that $2z \in \mbox{Im}
d^{1,v-2}$. Thus, our first step is to determine $\mbox{Im}
d^{1,v-2}$.

\noindent $\blacktriangle$ \emph{Matrix representation of
$d^{1,v-2}$:}

The labelling of the vertices and the edges of the graph is the the
same as described in the odd cycle case, as illustrated in Figure
\ref{odd cycle example}. The basis elements of $C^{1,v-2}(G)$ are
enhanced states $(s, c)$ where $s$  consists of one edge and
 $c$ is an assignment of $1$ or $x$ to the $v-1$ components of $[G:s]$.
 Since the total degree is $v-2$, all components are
assigned the value $x$ except one that is assigned the value $1$.
The basis element for which the present edge is $e_i$ and the
vertex that is assigned $1$ is $v_j$ is denoted by $b_i^j$. An
example is given in Figure \ref{even cycle example}.

\begin{figure}[h]
\centerline{\psfig{figure=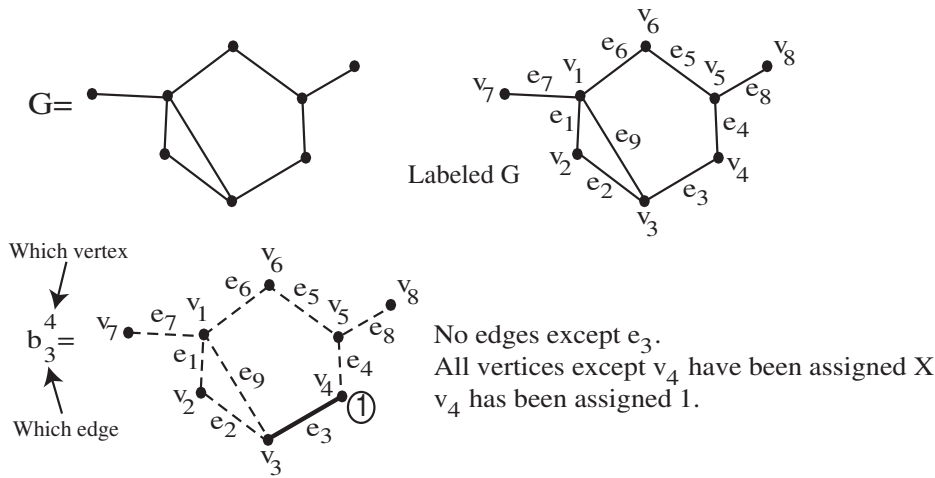,height=6.5cm}}
\caption{\textit{Notation for basis elements basis element in
 $C^{1,v-2}(G)$}} \label{even cycle example}
\end{figure}


We also need to describe the basis elements of the target space
$C^{2,v-2}(G)$. Let $(s,c)$ be an enhanced state in $C^{2,v-2}(G)$.
The graph $[G:s]$ has $v-2$ components, since adding two edges
automatically decreases the number of components by two (note that
$G$ contains no loop or multiple edges). Therefore for degree
reason, all the components are assigned the color $x$. The basis
element for which the present edges are $e_i$ and $e_j$, with $x$
assigned to all the components is denoted by $a_{ij}$ with $i<j$.

With these notations, we get a matrix representing $d^{1,v-2}$,
which is shown in the next paragraph.

\noindent $\blacktriangle$ We are now ready to exhibit an element
$z$ in $\ker d^2$ of degree $v-2$ such that $2z=0$  in $H^2(G)$ and
$z \neq 0$ in $H^2(G)$.


\noindent $\triangle$  We first prove that there exists an element
in $\mbox{Im} d^{1,v-2}$ with all coordinates even. This is the $2z$
we were looking for.

In the matrix representation of $d^{1,v-2}$, the coordinates of the
image of $b_i^j$ under the differential, in the basis
$(a_{ij})_{i,j}$ are listed in the $i^{\mbox{\tiny th}}$ row.

\bigskip

$ \hspace{-2cm}
\begin{array}{ccc}
\text{\emph{{\small Matrix of }}}d^{1,v-2}&
\begin{array}{ccc}
\begin{array}{c}
e_{1}\ \text{{\small and }} \\
\text{{\small another edge in cycle}}
\end{array}
&
\begin{array}{c}
e_{i}\text{, }e_{j}\text{ {\small in cycle, }} \\
i,j\neq 1
\end{array}
&
\begin{array}{c}
\text{{\small At least one edge}} \\
\text{{\small \ not in cycle}}
\end{array}
\\
\overbrace{
\begin{array}{ccccc}
a_{12} & a_{13} &  & \cdots  & a_{1,2s}
\end{array}
} & |\overbrace{
\begin{array}{ccc}
a_{23} & \cdots  & a_{2s-1,2s}
\end{array}
}| & \text{ {\small \ all other \ }}a_{ij}\text{ \ }
\end{array}
& \text{{\small Coeff}} \\
\begin{array}{c}
d(b_{1}^{1}) \\
d(b_{1}^{3}) \\
\\
\\
d(b_{1}^{2s}) \\
---- \\
d(b_{1}^{2s+1}) \\
\vdots  \\
d(b_{1}^{v}) \\
---- \\
d(b_{i}^{j}) \\
e_{i}\text{, }v_{j}\text{ in cycle} \\
i\neq 1 \\
---- \\
\\
\text{other }b_{i}^{j}
\end{array}
& \left[
\begin{array}{ccccccccccccc}
-1 &  &  &  & -1 & | &  &  &  & | &  &  &  \\
-1 & -1 &  &  &  & | &  &  &  & | &  &  &  \\
& -1 & \ddots  &  &  & | &  & 0 &  & | &  & M_{1} &  \\
&  & \ddots  & -1 &  & | &  &  &  & | &  &  &  \\
&  &  & -1 & -1 & | &  &  &  & | &  &  &  \\
- & - & - & - & - & -|- & - & - & - & -|- & - & - & - \\
&  &  &  &  & | &  &  &  & | &  &  &  \\
&  & 0 &  &  & | &  & 0 &  & | &  & M_{2} &  \\
&  &  &  &  & | &  &  &  & | &  &  &  \\
- & - & - & - & - & -|- & - & - & - & -|- & - & - & - \\
&  &  &  &  & | &  &  &  & | &  &  &  \\
&  & M_3 &  &  & | &  & M_4 &  & | &  & M_5 &  \\
&  &  &  &  & | &  &  &  & | &  &  &  \\
- & - & - & - & - & -|- & - & - & - & -|- & - & - & - \\
&  &  &  &  & | &  &  &  & | &  &  &  \\
&  & 0 &  &  & | &  & 0 &  & | &  & M_{6} &  \\
&  &  &  &  & | &  &  &  & | &  &  &
\end{array}
\right]  &
\begin{array}{c}
1 \\
1 \\
\vdots  \\
\vdots  \\
1 \\
-- \\
1 \\
\vdots  \\
1 \\
-- \\
0 \\
\vdots  \\
0 \\
-- \\
0 \\
\vdots  \\
0
\end{array}
\\
d(\sum_{1\leq j \leq v} b_1^j)= & \overline{
\begin{array}{ccccccccccccc}
-2 & -2 & \cdots  & \cdots  & -2 & \text{ \ }|\text{ \ } & \text{
}0 & \cdots  & 0 & \text{ \ }|\text{ \ } &
 & -2\epsi_{ij} &
\end{array}
}=2z &
\end{array}
$

where $\epsi_{ij}\in \{0,1\}$ for all $(i,j)\in J$ where $J$ is
the set of all $(i,j)$ such that at least one of the edges $e_i$,
$e_j$ is not in the cycle and $i<j$.

``Other $b_{ij}$" means either $e_i$ not in the cycle or $e_i$ is
in the cycle but $v_j$ isn't.

\bigskip
We need to explain why adding the rows of $M_1$ and $M_2$ always
yields a coordinate equal to $-2$ or $0$ on each $a_{ij}$, $(i,j)
\in J$. If $i=1$, each of these $a_{1j}$ is in the image of
exactly two $b_1^j$'s, the ones corresponding to a $1$ placed at
each endpoint of the edge $e_j$ (since each edge has two ends
under the assumption that there are no loops). If $i \neq 1$,
$a_{ij}$ is in the image of none of the $b_1^j$'s hence the
coordinate on the $a_{ij}$ with $i \neq 1$ are $0$.

It remains to explain the negative signs. Each $b_1^j$ is a basis
element coming from the state labelled $10....0$ (the present edge
is the first one in the ordering) so the label of the per-edge map
that adds the edge $i$ with $i \geq 2$ is $10..0*0..0$ with the
star in the $i^{\mbox{\tiny th}}$ position. The definition of the
differential says that when there is an odd number of $1$'s before
the star, the map is assigned a negative sign.

By adding the  first $v$  rows of this matrix, we see that all the
coordinates of $d\left( \sum_{j=1}^{v} b_1^{j}\right)$ are even so
we can call this element $2z$. Hence we have achieved our first goal
which is to find $2z \in \mbox{Im} d^{1,v-2}$.

Note that since $C^3(G)$ doesn't have torsion, this implies that $z
\in \ker d^{2,v-2}$.

\bigskip \noindent $\triangle$ It remains to show that  $z\notin \mbox{Im} d^{1,v-2}$.
Assume $z$ can be written $z=d(y)$ for some $y \in C^{1,v-2}(G)$. We
write the coordinates of $y$ in the same basis of $C^{1,v-2}(G)$ as
the one we previously used to write the matrix expression for
$d^{1,v-2}$. Namely, the basis for $C^{1,v-2}(G)$ we use can be
described as a partition $B_1 \sqcup B_2 \sqcup B_3$ where $B_1$ is
the set of basis elements for which the present edge is $e_1$, $B_2$
is the set of basis elements $b_i^j$ such that $e_i$ is an edge in
the cycle but is not $e_1$ and $v_j$ is a vertex in the cycle, and
$B_3$ is the set of all other basis elements $b_i^j$, i.e. either
the ones for which $e_i$ is an edge in the cycle but $v_j$ is not a
vertex in the cycle or $e_i$ is not in the cycle. Using this basis,
we can write $y$ as a linear combination of basis elements, labeling
its coordinates on elements of $B_1$ by $\alpha_i$, its coordinates
on elements of $B_2$ by $\beta_i$,and its coordinates on elements of
$B_3$ by $\gamma_i$, as illustrated in the following matrix
representation.



\bigskip

$\hspace{-2.5cm}
\begin{array}{ccc}
\text{\emph{{\small Matrix of }}}d^{1,v-2}&
\begin{array}{ccc}
\begin{array}{c}
e_{1}\ \text{{\small and }} \\
\text{{\small another edge in cycle}}
\end{array}
&
\begin{array}{c}
e_{i}\text{, }e_{j}\text{ {\small in cycle, }} \\
i,j\neq 1
\end{array}
&
\begin{array}{c}
\text{{\small At least one edge}} \\
\text{{\small \ not in cycle}}
\end{array}
\\
\overbrace{
\begin{array}{ccccc}
a_{12} & a_{13} &  & \cdots  & a_{1,2s}
\end{array}
} & |\overbrace{
\begin{array}{ccc}
a_{23} & \cdots  & a_{2s-1,2s}
\end{array}
}| & \text{ {\small \ all other \ }}a_{ij}\text{ \ }
\end{array}
& \text{{\small Coeff}} \\
\begin{array}{c}
d(b_{1}^{1}) \\
d(b_{1}^{3}) \\
\\
\\
d(b_{1}^{2s}) \\
---- \\
d(b_{1}^{2s+1}) \\
\vdots  \\
d(b_{1}^{v}) \\
---- \\
d(b_{i}^{j}) \\
e_{i}\text{, }v_{j}\text{ in cycle} \\
i\neq 1 \\
---- \\
\\
\text{other }b_{i}^{j}
\end{array}
& \left[
\begin{array}{ccccccccccccc}
-1 &  &  &  & -1 & | &  &  &  & | &  &  &  \\
-1 & -1 &  &  &  & | &  &  &  & | &  &  &  \\
& -1 & \ddots  &  &  & | &  & 0 &  & | &  & M_{1} &  \\
&  & \ddots  & -1 &  & | &  &  &  & | &  &  &  \\
&  &  & -1 & -1 & | &  &  &  & | &  &  &  \\
- & - & - & - & - & -|- & - & - & - & -|- & - & - & - \\
&  &  &  &  & | &  &  &  & | &  &  &  \\
&  & 0 &  &  & | &  & 0 &  & | &  & M_{2} &  \\
&  &  &  &  & | &  &  &  & | &  &  &  \\
- & - & - & - & - & -|- & - & - & - & -|- & - & - & - \\
&  &  &  &  & | &  &  &  & | &  &  &  \\
&  & M_3 &  &  & | &  & M_4 &  & | &  & M_5 &  \\
&  &  &  &  & | &  &  &  & | &  &  &  \\
- & - & - & - & - & -|- & - & - & - & -|- & - & - & - \\
&  &  &  &  & | &  &  &  & | &  &  &  \\
&  & 0 &  &  & | &  & 0 &  & | &  & M_{6} &  \\
&  &  &  &  & | &  &  &  & | &  &  &
\end{array}
\right]  &
\begin{array}{c}
\alpha_1 \\
\alpha_3 \\
\vdots  \\
\alpha_{2s-1} \\
\alpha_{2s} \\
-- \\
\alpha_{2s+1}  \\
\vdots  \\
\alpha_{v} \\
-- \\
\beta_1 \\
\vdots  \\
\beta_m \\
-- \\
\gamma_1 \\
\vdots  \\
\gamma_{m'}
\end{array}
\\
d(y)= & \overline{
\begin{array}{ccccccccccccc}
-1 & -1 & \cdots  & \cdots  & -1 & \text{ \ }|\text{ \ } & \text{
}0 & \cdots  & 0 & \text{ \ }|\text{ \ } &
 & -\delta_{ij} &
\end{array}
}=z &
\end{array}
$

Note that there is no $b_1^2$ in $B_1$ since $b_1^2$ and $b_1^1$
would be the same (and no corresponding $\alpha_2$ coefficient).
Hence the first block matrix with the $-1$'s is a
$(2s-1,2s-1)$-matrix.

$\bigskip $

If we now read this by columns, what we get for the two first
blocks of columns in the matrix representation, i.e. the columns
corresponding to basis elements with both edges in the cycle, is
the following:

$\hspace{-2.5cm}
\begin{array}{c}
\text{ \ }
\begin{array}{c}
\text{ }
\begin{array}{c}
\text{Column 1}(a_{12})\text{:} \\
\text{Column 2}(a_{13}) \text{:}\\
\vdots  \\
\\
\text{Column 2s-1}(a_{1,2s})\text{:}
\end{array}
\\
\text{ }
\begin{array}{c}
------- \\
\text{Column 2s} (a_{23})\text{:}\\
\vdots  \\
\\
\text{Column }\left(\tiny{
\begin{array}{c}
2s \\
2
\end{array}}
\right) (a_{2s-1,2s})\text{:}
\end{array}
\end{array}
\text{\ \ }
\begin{array}{c}
\begin{array}{cccccccc}
-\alpha _{1} & -\alpha _{3} &  &  &  & | & +S_{1} & =-1\text{ \ \ \ } \\
& -\alpha _{3} & -\alpha _{4} &  &  & | & +S_{3} & =-1\text{ \ \ \ } \\
&  & \ddots  & \ddots  &  & | &  & \vdots  \\
&  &  & -\alpha _{2s-1} & -\alpha _{2s} & \text{ \ }|\text{ \ } &
+S_{2s-2}
& =-1\text{ \ \ \ } \\
-\alpha _{1} &  &  &  & -\alpha _{2s} & | & +S_{2s-1} & =-1\text{
\ \ \ }
\end{array}
\\
\text{ \ \ \ }
\begin{array}{cccccccc}
-- & -- & -- & -- & -- & -|- & -- & \text{ \ \ \ } \\
\text{ \ \ \ \ \ } & \text{ \ \ \ \ \ } & \text{ \ \ \ \ \ } &
\text{ \ \ \
\ \ \ \ } & \text{ \ \ \ \ \ } & | & +S_{2s} & =0\text{ \ \ \ } \\
&  & 0 &  &  & | &  & \vdots  \\
&  &  &  &  & | &  & \text{ \ \ \ } \\
&  &  &  &  & | & +S_{\left( \tiny{\begin{array}{c}
2s \\
2
\end{array}}
\right) } & =0\text{ \ \ \ }
\end{array}
\end{array}
\\
\text{ \ \ \ \ \ \ \ \ \ \ \ \ \ \ \ }
\begin{array}{c}
\text{\ \ \ \ \ \ \ \ \ \ \ }
\end{array}
\overline{
\begin{array}{cccccccc}
-2(\alpha _{1} & +\alpha _{3} & +\alpha _{4} & \cdots  & \cdots  &
+\alpha _{2s})| & +S_{\beta } & \text{\ \ \ \ \ }=2s-1
\end{array}
}
\end{array}
$

where the $S_{i}$'s the result of the multiplication of $M_3$ and
$M_4$ by the $\beta _{i}$'s, read by columns. Their sum $S_{\beta
}$ is a linear combination of $\beta _{i}$'s with coefficients in
$\mathbb{Z}$. It suffices to prove that these coefficients are all
even to get a contradiction, since this will imply that the left
hand side is even while the right hand side is odd. This will be
achieved by showing that there are exactly two non-zero entries in
each row of the matrix $M=\left[
\begin{array}{ccc}
M_{3} & | & M_{4}
\end{array}
\right]$, and that these entries are $\pm 1$.

Indeed, for any $b_i^j$ with $e_i, v_j$ in the cycle, there are
exactly two ways to add an edge adjacent to the component with the
color $1$ under the condition that this edge is in the cycle. The
coordinates of $d(b_i^j)$ on these two basis elements is $\pm 1$
and appear in $M$ . The other basis elements in the image of
$b_i^j$ under the differential will be have at least one edge not
in the cycle so the corresponding coordinate will appear in $M_5$.
\end{proof}

\subsection{Connection with
 Khovanov cohomology for links}

\begin{figure}[h]
\begin{center}
\scalebox{.6 }{\includegraphics{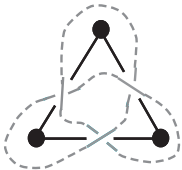}}
\end{center}
\caption{\textit{The left handed trefoil knot and its Tait graph
$P_3$}} \label{TaitTrefoilP3}
\end{figure}

Recall that there is a bijection, originally described by Tait,
between the set of unoriented alternating link diagrams and the set
of plane graphs. Indeed, one can draw a crossing on each edge of a
plane graph according to the following diagram
$\vcenter{\hbox{\epsfig{file=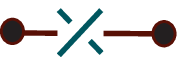}}}$ \footnote{ Our
convention of creating a link diagram out of a graph corresponds to
negative graph (signed graph with all negative edges). However, we
do not deal in this paper with signed graphs so our convention
should not lead to confusion. Our choice is dictated by the fact
that we want a $0$-resolution of the crossing in the link case to
correspond to the case when the edge is absent in the graph case.},
and then connect along the edges, as illustrated by the dotted lines
in Figure \ref{TaitTrefoilP3}. This way we obtain an unoriented
alternating link diagram.

Conversely, we can obtain a plane graph from the diagram of an
alternating link using the chess-board coloring according to the
convention $\vcenter{\hbox{\epsfig{file=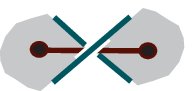}}}$ (for
details see \cite{P07}).

The version of Khovanov cohomology that we are going to use here
is the one for framed unoriented links defined by Viro in Section
6 of \cite{V04}. It categorifies the augmented Kauffman bracket
polynomial, i.e. the version of the Kauffman bracket polynomial
which takes the value $1$ on the empty link. Since different
notations have been used to describe this categorification, we are
first going to specify which ones we use. The crossing are
labelled $1$ to $n$ and a $0$-marker (resp. a $1$-marker)
associated to a crossing indicates that we perform an
$A$-smoothing (resp. an $A^{-1}$-smoothing) at this crossing. An
$A$-smoothing is what Khovanov calls a 0-resolution and what
Bar-Natan calls a 0-smoothing. They correspond to the positive
markers of Viro.

\begin{figure}[h]
\begin{center}
\scalebox{.65 }{\includegraphics{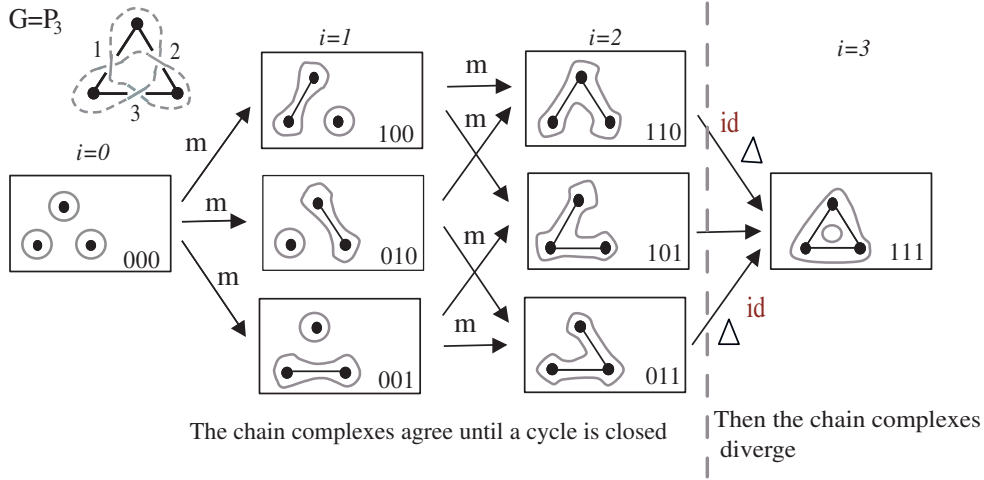}}
\end{center}
\caption{\textit{Comparison of the chain complexes for the left
handed trefoil knot and its Tait graph $P_3$ when the algebra is
$\Atwo$}} \label{RelKhovLinks}
\end{figure}

Figure \ref{RelKhovLinks} illustrates the relation between the chain
complexes. In this figure, we adopt a presentation of our chain
complex similar to the one in \cite {BN02}.

\begin{theorem}\label{Let $D$ be the diagram}
Let $D$ be the diagram of an unoriented framed alternating link
and let $G$ be its Tait graph. Let $\ell$ be the length of the
shortest cycle in $G$. For all $i<\ell-1$, we have

\[
H^{i,j}_{Graph}(G) \cong H_{p,q}(D)
\]
with $\left \{\begin{array}{l}
p=V(G)-i-2j,\\
q=E(G)-2V(G)+4j.
\end{array}
\right.$

where $H_{p,q}$ is the homology groups of the version of Khovanov
cohomology theory for unoriented framed link defined in
\cite{V04}.\footnote{If we use the previous definition in
\cite{V02}, we get 
$H^{ij}_{Graph}(G) \cong H_{a,b}(D)$
 with $\left \{\begin{array}{l}
a=E(G)-2i,\\
b=E(G)-2V(G)+4j.
\end{array}
\right.$ In both Viro's definitions the degree is the same $q=b$
and corresponds to the $A$-degree of the Kauffman bracket
polynomial.}

Furthermore $Tor(H^{i,j}(G) = Tor(H_{p,q}(D))$ for $i=\ell-1$.
\end{theorem}

\begin{proof}
For each state $\ga \in \{0,1\}^n$, as long as adding edges does not
close a cycle in the graph, the connected components of the graph
and the link case naturally correspond to each other. See Figure
\ref{RelKhovLinks} for an example. Viro doesn't mention the
$\Z$-algebra $\Atwo$ in his construction but he assigns $+$ and $-$
signs to the circles and then combines them when two discs merge
using the obvious multiplication: $+*+=+, +*-=-*+=-$ and $-*-=+$.
With this multiplication, the free $\Z$-module with basis + and - is
a $\Z$-algebra isomorphic to $\Atwo$ under the correspondence $ +
\leftrightarrow 1, - \leftrightarrow x$. In both cases we build the
chain groups by assigning a copy of $\A_2=\Atwo$  to each component,
taking tensor products over components and direct sum along columns.
Therefore the chain groups are isomorphic as long as the spanning
subgraphs are trees.

The differentials also agree: As long as adding edges does not
close a cycle in the graph, changing one marker from 0 to 1
corresponds to connecting discs in the links case and to
connecting components in the graph case so in both cases the new
labels are computed using the multiplication of $\A_2$. This is
enough to see that the cohomology groups are isomorphic for
$i<\ell-1$.

The definition for the heights and gradings are different in these
two chain complexes so we need express the relation between them.
Let $i$ be an integer with $i<\ell-1$. Recall that an enhanced
state $S$ in the link cohomology is a smoothing of the link
diagram according to markers with each of the resulting circles
colored by $1$ or $x$. Given an enhanced state $S$ in the link
cohomology $p(S)$ is defined to be the number of $1$ minus the
number of $x$ \cite{V04}. This can be written as $p(S)=\#1
-\#x=\#1 +\#x-2\#x=\#$circles $-2j$ since in the graph cohomology
the grading $j$ is the number of $x$. Also, note that since
$i<\ell-1$, the number of circles in the link cohomology and the
number of connected components of the graph in the corresponding
graph cohomology are equal. Therefore, $p(S)=\#
\mbox{components}([G:s]) - 2j$. Since adding edges didn't close
any cycle by our assumption on $i$, the number of connected
components of the graph is equal to the number of vertices $V(G)$
minus the number of edges $i$. This proves that $p(S)=V(G)-i- 2j$.

We now deal with $q(S)$. Given an enhanced state $S$ in the link
cohomology  $q(S)$ is defined to be the number of $0$-resolution
of a crossing minus the number of $1$-resolution of a crossing
minus $2p(S)$ \cite{V04}. This can be written
$q(S)=\#\mbox{0-resolution}- \#\mbox{1-resolution}
-2p(S)=\#\mbox{0-resolution}+
\#\mbox{1-resolution}-2\#\mbox{1-resolution}
-2p(S)=\#\mbox{crossings} -2(\#\mbox{1-resolution})-2p(S)$ since
the number of crossings is equal to the number of edges of $G$ and
the number of $1$-resolutions is equal to $i$, the number of edges
in $S$. Therefore $q(S)=E(G)-2i- 2(V(G)-i- 2j)$. This proves that
$q(S)=E(G)-2V(G)+4j$.

Now, note that we can get a slightly better result for torsion
than for the whole homology groups. Indeed, if we want to compute
only the torsion part of some homology groups, we look only at the
image, not the kernel of the differentials at that height and the
images in chain complexes agree until $i=\ell -1$.
\end{proof}

\section{Computations for algebras $\A_m$ and their deformations}\
\label{section the algebra is Am}


One can conjecture that results similar to that of Section
\ref{section algebra is A2} (e.g. Theorem \ref{torsion theorem
A2}) hold generally for algebras $\A_m$. In particular

\begin{conjecture} \label{torsion conjecture Am} The cohomology $H_{A_m}^{**}(G)$ of a
graph $G$ contains a torsion part if and only if  $G$ has no loops
and contains a cycle of order greater than or equal to $3$. In
this case, $H^{**}_{\A_m}(G)$ has a torsion of order dividing $m$.
\end{conjecture}

In this section we show several examples confirming the conjecture
and in particular we compute $H^{i,j}_{\A_m}(P_3)$.

As we noted in the previous section, graph cohomology for $\A_2$ can
be used to approximate classical $sl(2)$ Khovanov cohomology. We can
speculate that our more general calculation can be used to
approximate Khovanov-Rozansky $sl(m)$ homology \cite{KR04} (in that
case $\A=\A_m=\Z[x]/x^m)$ or its deformation corresponding to
algebra $\A=\A_{p(x)} = \Z[x]/(p(x) $, where $p(x)$ is a polynomial
in $x$ (compare \cite{Gor}). Almost no computation is completed for
Khovanov-Rozansky $sl(m)$ homology and speculations concern mostly
the free part of them (compare \cite{GuSchVa,DGR}). Our calculations
can be very usefull for predicting and verifying future results.

\subsection{$H^{i,j}_{\A_m}$ for a triangle}

We compute here the graph cohomology $H^{i,j}_{\A_m}(P_3)$,
including the case of the polynomial ring $\A_{\infty}=\Z[x]$. We
formulate our main theorem in two parts describing torsion part
and free part
 separately.

Recall that the \emph{Poincar\'{e} series} of a sequence of graded
$\Z$-modules $\{M^i=\oplus_j M^{i,j}\}_i$ where $M^{i,j}$ is the set
of homogeneous elements of $M^i$ of degree $j$, is defined to be
two-variable series $\sum_{i,j} t^i q^j \dim_{\mathbb{Q}}(M^{i,j}
\otimes_{\Z} \mathbb{Q})$. With our definition of the graded
dimension, this can be rewritten $\sum_{i} t^i qdim(M^i)$. The
Poincar\'{e} series of a sequence of graded $\Z$-modules is a
convenient way to store the free ranks of the $M^{i,j}$'s.

Note that if two sequence $\{M^i\}_i$ and $\{M'^i\}_i$ of
\emph{free} graded
 $\Z$-modules have the same Poincar\'{e} series, then
the groups $M^{i,j}$ and $M'^{i,j}$ are isomorphic
for all $i,j$.

\begin{theorem} When the algebra is $\A_m$, the
cohomology groups of $P_3$ are completely determined by:\ \
\label{thm Am P3}
\begin{enumerate}
\item[(Torsion)] $Tor(H^{*,*}_{{\A}_m}(P_3))=
H^{1,m}_{{\A}_m}(P_3)= \Z_m$. \item[(Free)] The Poincar\'{e}
polynomial of $H^{*,*}_{{\A}_m}(P_3)$ is equal to
$(q+q^2+...+q^{m-1})^3 + t(q+q^2+...+q^{m-1})$.
\end{enumerate}
\end{theorem}


\begin{theorem} When the algebra is $\A_{\infty}$, the
cohomology groups of $P_3$ are completely determined by: \
\label{thm A infty}
\begin{enumerate}
\item[(Height 0)] $H^{0,*}_{\A_{\infty}}(P_3)$ is the free abelian group
$(\A')^{\otimes 3}$, where $\A'=
\Z\{1\} \oplus \Z\{2\} \oplus \Z\{3\} \oplus
\cdots=\oplus_{j=1}^{\infty}\Z\{j\}$. \item[(Height 1)]
$H^{1}_{\A_{\infty}}(P_3)=\A'$ i.e. $H^{1,0}_{\A_{\infty}}(P_3)= 0$
and $H^{1,j}_{\A_{\infty}}(P_3)= \Z$ for any $j>0$.
\end{enumerate}
\end{theorem}
By Remark \ref{Remark6}, $H^{i}_{\A_{\infty}}(P_3)=0 $ for all
$i\geq 2$. Thus, tatements (0) and (1) together are equivalent to
saying that the cohomology groups $H^{*,*}_{\A_{\infty}}(P_3)$ are
free abelian with Poincar\'{e} series $(q+q^2+q^3+...)^3 +
t(q+q^2+q^3+...)$.


It is not clear for which $\A$ the group
 $H^{0,*}_{\A}(P_3)$ is isomorphic to $(\A')^{\otimes 3}$ and
$H^{1,*}_{\A}(P_3)$ is isomorphic to $\A'$.

\begin{corollary}
 Let $P_v$ be an $v$-gon with $v \geq 3$; then
\begin{enumerate}
\item[(Torsion)] $Tor(H^{v-2,*}_{{\A}_m}(P_v))=
H^{v-2,m}_{{\A}_m}(P_v)= \Z_m$. \item[(Free)] The Poincar\'{e}
polynomial of $H^{v-2,*}_{{\A}_m}(P_v)$ is equal to
$t^{v-2}(q+q^2+...+q^{m-1})$.
\end{enumerate}
\end{corollary}

\begin{proof}
 By corollary
\ref{Polygon in terms of H1},
$H^{v-2}_{{\A}_m}(P_v)=H^{1}_{{\A}_m}(P_3)$.
\end{proof}

We prove Theorems \ref{thm Am P3} and \ref{thm A infty}
simultaneously. The idea of the proof is to analyze carefully the
chain maps $0 \to C^{0,j} \stackrel{d^{0,j}}{\to} C^{1,j}
\stackrel{d^{1,j}}{\to} C^{2,j} \stackrel{d^{2,j}}{\to} C^{3,j} \to
0$ for $P_3$ and ${\A}_{\infty}$. Recall that $d^{i,j}$ denotes the
restriction of $d^i$ to the elements of degree $j$ of $C^i$. Precise
knowledge of ${d^0}$ and ${d^1}$ will allow us also to obtain result
for ${\A}_{m}$.

\begin{figure}[h]
\begin{center}
\scalebox{.9 }{\includegraphics{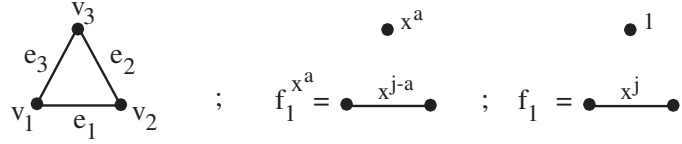}}
\end{center}
\caption{\textit{Basis elements of
$C^{1,j}_{{\A}_{\infty}}(P_3)$.}} \label{basisC1j}
\end{figure}

 We first describe the basis elements of
$C^{1,j}_{{\A}_{\infty}}(P_3) = \Z^{3(j+1)}$. Each basis element is
made of two connected components, the one with the edge and the one
consisting of the remaining vertex. Each component is colored by a
power of $x$ such that the sum of the exponents is $j$ since we are
dealing with degree $j$ elements. We denote by $f_i^{x^a}$ the basis
element for which the present edge is $e_i$ and the color on the
isolated vertex is $x^a$. The color on the other component does not
need to be specified in the notation since it is  forced to be
$x^{j-a}$. To alleviate the notation we denote $f_i^1$ simply by
$f_i$ for $i=1,2,3$. Examples are shown in Figure \ref{basisC1j}.


\begin{lemma}[Main Lemma]
\label{main lemma} Consider  $C_{{\A}_{\infty}}^{*,j}(P_3)$ with
$j>0$. Then
\begin{enumerate}
\item [(a)]$f_1+f_2+f_3$, $f_1^{x^a}-f_1$, $f_2^{x^a}-f_2$ and
$f_3^{x^a}-f_3$ ($a\geq 1$) form a basis of $ker(d^{1,j})$.

\item [(b)]Let $R$ denote the subgroup of
$C_{\A_{\infty}}^{1,j}$ generated by the elements $f_1+f_2+f_3 +
j(f_1^x-f_1)$ and all elements of the form $(f_i^{x^a}-f_i) -
a(f_1^x-f_1)$ with $a \geq 1$, $i=1,2,3$. Then $\mathfrak{R} \subset
Im(d^{0,j})$.

\item[(c)] Elements from (b) form a basis of $Im(d^{0,j})$ and in
particular, $H^{1,j}_{{\A}_{\infty}}(P_3)=
ker(d^{1,j})/Im(d^{0,j})$ is isomorphic to $\Z$ and generated
by the class of $f_1^x-f_1$.\\
 Furthermore $f_1+f_2+f_3 \equiv
-j(f_1^x-f_1)$ in $H^{1,j}_{{\A}_{\infty}}(P_3)$.
\end{enumerate}
\end{lemma}

\begin{proof}[Proof of lemma \ref{main lemma}.]
Let us first determine $C_{{\A}_{\infty}}^{0,j}(P_3)$,
 and $C_{{\A}_{\infty}}^{2,j}(P_3)$
and introduce notations for their basis elements.

The vertices of $P_3$ are labelled as in Figure \ref{basisC1j}.
The triplet $(x^a, x^b, x^c)$ obviously denotes the basis element
of $C_{{\A}_{\infty}}^{0,j}(P_3)$ where the color $x^a$ is
associated to the first vertex, $x^b$ is associated to the second
vertex and $x^c$ is associated to the third vertex. It is assumed
whenever we use this notation that $a+b+c=j$.

The basis elements of $C_{{\A}_{\infty}}^{2,j}(P_3)$ have two
edges and only one connected component. Therefore, the color
assigned to this component is $x^j$. If the two edges are $e_i$
and $e_k$, the corresponding basis element is denoted
$(x^j)(e_i,e_k)$. It is clear that $C^{2,j}(P_3)=\Z^3$ with basis
$(x^j)(e_1,e_2)$, $(x^j)(e_2,e_3)$ and $(x^j)(e_1,e_3)$

 \bigskip

\emph{Proof of part (a) of the lemma.} The elements $f_i^{x^a}$ with
$i=1,2,3$ and $ 0\leq a \leq j$ form a basis of
$C^{1,j}=\Z^{3(j+1)}$.  We will also use the more convenient basis
formed by the elements $f_2$, $f_3$, $f_1+f_2+f_3$, $f_1^{x^a}-f_1$,
$f_2^{x^a}-f_2$ and $f_3^{x^a}-f_3$ ($a>0$). The map $d^{2,j}$ is an
epimorphism with $\ker (d^{2,j})=\Z^2\{j\}$ generated by
$(q^j)(e_1,e_2)-(q^j)(e_2,e_3)$ and $(q^j)(e_1,e_2) +
(q^j)(e_1,e_3)$. The map $d^{1,j}$ is an epimorphism onto
$\ker(d^{2,j})$ and in fact $d^{1,j}(f_2)=
(x^j)(e_1,e_2)-(x^j)(e_2,e_3)$ and $d^{1,j}(f_3) = (x^j)(e_1,e_2) +
(x^j)(e_1,e_3)$ span $\ker(d^{2,j})$. Therefore the $3j+1$ elements
$f_1+f_2+f_3$, $f_1^{x^a}-f_1$, $f_2^{x^a}-f_2$ and $f_3^{x^a}-f_3$
($a>0$) form a basis of $\ker (d^{1,j})$.

Note that the fact that $d^{1,j}$ is surjective shows that
$H_{{\A}_{\infty}}^{2,j}(P_3)=0$, a fact that we could also get
using Theorem \ref{vanishing theorem}.

 \bigskip
 \emph{Proof of  part (b) of the lemma.} We need to
analyze carefully the image of $d^{0,j}$. To simplify our
considerations we will use $\tau$, the automorphism of $C^{i,j}$
induced by the $2\pi/3$ counter-clockwise rotation by $2\pi/3$.
$d^{0,j}$ commutes with $\tau$, i.e. $\tau d^{0,j} = d^{0,j}
\tau$. In fact action by the full dihedral symmetry group of the
triangle commutes with $d^{0,j}$. We also have $\tau^2(f_1^{x^a})=
\tau(f_2^{x^a})=f_3^{x^a}$. An example is shown in Figure
\ref{basisC1jTau}.
\begin{figure}[h]
\begin{center}
\scalebox{.9 }{\includegraphics{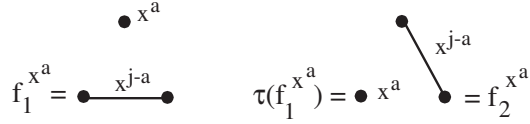}}
\end{center}
\caption{\textit{Image of the basis elements of
$C^{1,j}_{{\A}_{\infty}}(P_3)$
 under $\tau$.}}
\label{basisC1jTau}
\end{figure}

 With this in mind, consider Im$d^{0,j}_{R1}$ which
denotes the image under $d^{0,j}$ of the subgroup  $C^{0,j}(R1)$
of $C^{0,j}$ generated by basic elements
 $(x^{a_1},x^{a_2},x^{a_3})$ with at least one $a_i$ equal to
$0$. Consider the generic element $(1,x^a,x^b)$ (all other are
obtained from this one by rotation and taking different $0\leq a
\leq j$ with $b=j-a$). We have $d^{0,j}((1,x^b,x^a))= f_1^{x^a} +
f_2 +f_3^{x^b}$ (1). This relation and its image under $\tau^2$
with the role of $a$ and $b$ reversed, that is $d^0((x^a,x^b,1))=
f_1 + f_2^{x^a} +f_3^{x^b}$ (2), are illustrated in Figure
\ref{relationsR1}.


\begin{figure}[h]
\begin{center}
\scalebox{.45 }{\includegraphics{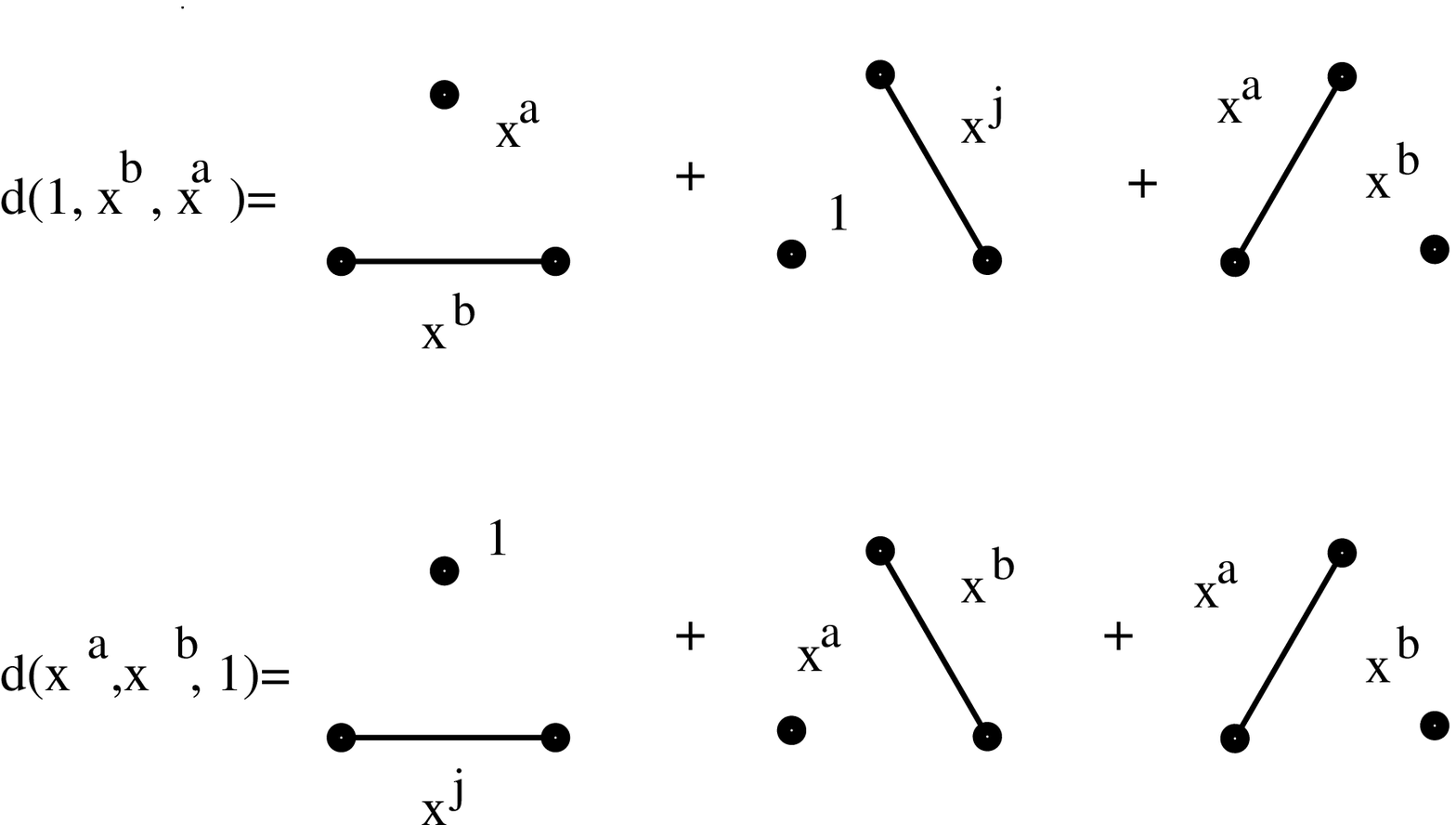}}
\caption{\vtop{\hbox{$ d^0((1,x^b,x^a))= f_1^{x^a} + f_2 +f_3^{x^b}$
\hspace{5mm} (1)}\hbox{$d^0((x^a,x^b,1))= f_1 + f_2^{x^a}+f_3^{x^b}$
\hspace{5mm} (2).}} \label{relationsR1}}
\end{center}
\end{figure}


Subtracting relation (2) from relation (1) shows that the elements
of the form $(f_1^{x^a}-f_1) - (f_2^{x^a}-f_2)$ are in
Im$d^{0,j}_{R1}$ for any $a$. Applying $\tau$ to this relation, we
get that the elements of the form $(f_i^{x^a}-f_i) -
(f_{i+1}^{x^a}-f_{i+1})$ are in Im$d^{0,j}_{R1}$ for any $a$ and any
$i \mod 3$. We call this (3). We write $(f_i^{x^a}-f_i)
\stackunder{R1}{\equiv} (f_{i+1}^{x^a}-f_{i+1})$ to express the
  fact that these terms are equal in
  $\ker(d^{1,j})/Imd_{R1}^{0,j}$.

Now, note that $f_1^{x^a} + f_2 +f_3^{x^b} \stackunder{R1}{\equiv}
0 $ by (1). The left hand side is equal to $f_1^{x^a}-f_1+f_1 +
f_2 +f_3^{x^b}-f_3+f_3$, therefore $(f_1^{x^a}-f_1) +
(f_3^{x^b}-f_3)+f_1 + f_2+f_3\stackunder{R1}{\equiv} 0$. Using (3)
with $i=3$, we get $(f_1^{x^a}-f_1) + (f_1^{x^b}-f_1)+f_1 +
f_2+f_3\stackunder{R1}{\equiv} 0$.

From above calculations, one can see that
\[
\begin{array}{ll}
(f_i^{x^a}-f_i) - (f_{1}^{x^a}-f_{1})\in  Imd_{R1}^{0,j} ~  \mbox{
for any }a \mbox{ and } i=1, 2, 3. &
 \hspace{3mm} (3') \\
(f_i^{x^a}-f_i) + (f_i^{x^b}-f_i) + f_1 +f_2 +f_3\in Imd_{R1}^{0,j}
~ \mbox{  for any } i \mbox{ and } a+b=j& \hspace{3mm} (4)
\end{array} \]

 Consequently, $\ker(d^{1,j})/$Im$d_{R1}^{0,j}$ is generated by
the classes of $f_1 +f_2 +f_3$ and of $f_1^{x^a}-f_1$ for $a \leq
\frac{j}{2}$.

\bigskip
Now consider the subgroup $C^{0,j}(R2)$ of $C^{0,j}$ generated by
enhanced states $(x^{a-1},x^b,x)$ for $1 \leq a-1 \leq j-1$ with
$a+b=j$. We have $d^{0,j}(x^{a-1},x^b,x) = f_1^{x} + f_2^{x^{a-1}}
+f_3^{x^b}$. Let Im$d^{0,j}_{R1 \cup R2}$ denote the image under
$d^{0,j}$ of the subgroup  $C^{0,j}(R1 \cup R2)$ of $C^{0,j}$
generated by enhanced states $(x^{a_1},x^{a_2},x^{a_3})$ with at
least one $a_i$ equal to $0$ or to $1$. We rewrite the previous
relation  $f_1^{x} + f_2^{x^{a-1}} +f_3^{x^b} \equiv_{R1 \cup R2}
0$ to mean that the left hand side and the right hand side are
equal in $\ker(d^{1,j})/Imd_{R1\cup R2}^{0,j}$. Adding and
substracting $f_1$, $f_2$ and $f_3$ to the left hand side leaves
it unchanged so $(f_1^{x}-f_1)+f_1 + (f_2^{x^{a-1}}-f_2)+f_2
+(f_3^{x^b}-f_3)+f_3 \stackunder{R1\cup R2}{\equiv} 0$. Applying
(3') to the terms with indices 2 and 3 yields $(f_1^{x}-f_1)+
(f_1^{x^{a-1}}-f_1) +(f_1^{x^b}-f_1)+f_1 +f_2+f_3
\stackunder{R1\cup R2}{\equiv} 0$. We now add and substract
$(f_1^{x^{a}}-f_1)$ and then apply (4) since $a+b=j$. The result
is $(f_1^{x}-f_1)+ (f_1^{x^{a-1}}-f_1) -(f_1^{x^{a}}-f_1)
\stackunder{R1\cup R2}{\equiv} 0$ i.e. $(f_1^{x^{a}}-f_1)
\stackunder{R1\cup R2}{\equiv} (f_1^{x}-f_1)+
(f_1^{x^{a-1}}-f_1)$.

Applying this several times we get:

\[f_1^{x^{a}}-f_1 \stackunder{R1\cup R2}{\equiv}
a(f_1^{x} -f_1) \hspace{7mm} (5) \] and in particular
\[f_1^{x^{j}}-f_1 \stackunder{R1\cup R2}{\equiv} j(f_1^{x} -f_1) \hspace{7mm} (6)\]
Notice that for any $a$, applying (5) to the relation
$(f_1^{x^a}-f_1) + (f_1^{x^b}-f_1) \stackunder{R1\cup R2}{\equiv}
-(f_1 +f_2 +f_3)$ yields
\[j(f_1^{x} -f_1) \stackunder{R1\cup R2}{\equiv} -(f_1 +f_2 +f_3) \hspace{7mm} (7) \]
This completes the proof of part (b) of the lemma since it
consists of (7) and (3') substituted in (5).

\bigskip
\emph{Proof of  part (c) of the lemma.} We just proved that
${R}\subset $Im$d^{0,j}_{R1 \cup R2} \subset $Im$d^{0,j}$, and we
have, from (a) that $\ker(d^{1,j})/R \cong \Z$ is generated by the
class of $(f_1^{x}
-f_1)$.\\
To complete our proof of Main Lemma we need to show that taking the
quotient by Im$d^{0,j}$ rather than $R$ doesn't introduce new
relations, that is Im$d^{0,j}=R$. It is enough to observe that the
images of all the basis elements $d^{0,j}(x^{a_1},x^{a_1},x^{a_1})$,
with $a_1 +a_2 +a_3=j$,
 are in $R$. In fact, in
$\ker(d^{1,j})/R$ we have: $d^{0,j}((x^{a_1},x^{a_2},x^{a_3}))
 = f_1^{x^{a_3}} + f_2^{x^{a_1}}
+ f_3^{x^{a_2}} \stackunder{R}{\equiv}  (f_1^{x^{a_3}}- f_1) +f_1 +
(f_1^{x^{a_1}}-f_1) +f_2 + (f_1^{x^{a_2}}-f_1) +f_3
\stackunder{R}{\equiv} a_3(f_1^x-f_1) + a_1(f_1^x-f_1) +
a_2(f_1^x-f_1) + (f_1+f_2+f_3) \stackunder{R}{\equiv} j(f_1^x-f_1) +
(f_1+f_2+f_3) \stackunder{R}{\equiv} 0$ by (7).
\end{proof}

 We are now ready to prove Theorems
\ref{thm Am P3} and \ref{thm A infty}. Immediately from Lemma
\ref{main lemma} we get for $j>0$ that
$H^{1,j}_{{\A}_{\infty}}(P_3)=\Z$ which combined with
$H^{1,0}_{{\A}_{\infty}}(P_3)=0$ shows that
$H^{1,*}_{{\A}_{\infty}}(P_3)$ is a free abelian group with
Poincar\'{e} series $t(q+q^2+...)$.

$H^{0,*}_{{\A}_{\infty}}(P_3)$ can be deduced directly from
$H^{1,*}_{{\A}_{\infty}}(P_3)$ by using the formula (\ref{chromatic
as graded Euler Char}) in Section \ref{section graph cohomology}
that expresses the chromatic polynomial as the graded Euler
characteristic of the chain complex, namely $\sum _{0\leq i\leq
n}(-1)^{i} q\dim (H^{i})=P_G(q\dim \A)$. In our specific case, this
becomes $q\dim (H^{0}(P_3))-q\dim
(H^{1}(P_3))=P_{P_3}(1+q+q^2+...)$. From this we can derive that the
Poincar\'{e} series of $H^{0,*}_{{\A}_{\infty}}(P_3)$
 equal to $(q+q^2+...)^3$.


 Thus the Poincar\'{e} series of
$H^{*,*}_{{\A}_{\infty}}(P_3)$ is equal to $(q+q^2+...)^3 +
t(q+q^2+...)$. This proves Theorem \ref{thm A infty}.

\bigskip
In order to prove Theorem \ref{thm Am P3}, we deduce  the case
$\A=\A_{m}$ from the case $\A=\A_{\infty}$ by adding the relation
$x^m=0$ (see Theorem 36 for more general setting).

 For $j<m$ we have immediately
$H^{*,j}_{{\A}_{m}}(P_3)= H^{*,j}_{{\A}_{\infty}}(P_3)$.

For $j=m$ we have $f_1=f_2=f_3=f_1^{x^j}=f_2^{x^j}=f_3^{x^j}=0$,
otherwise the proof of Main lemma works without change. We get as
before that the class of $f_1^x$ generates $\ker(d^1)/(Im(d^0)$.
Furthermore from $f_1+f_2+f_3=0$ follows that $jf_1^x=0$ and
$H^{1,j}_{{\A}_{m}}(P_3)=\Z_j$.

For $j>m$ we obtain immediately that $\ker(d^{1,j})\subset
Im(d^{0,j})$. Therefore $H^{1,j}_{{\A}_{m}}(P_3)=0$. Theorem
\ref{thm Am P3} follows.

\subsection{More calculations}
For a general polygon $P_v$, $v\geq 3$ one can compute
$H^{*,*}_{{\A}_{m}}$ provided that we know
$H^{1,*}_{{\A}_{m}}(P_{v'})$ for $v'\leq v$, as explained in
Corollary \ref{Polygon in terms of H1}. The result in the classical
case of $m=2$ is stated in Theorem \ref{Homology polygon A2} and
could also be deduced from Khovanov result on cohomology of
$(2,v)$-torus links \cite{K00}. We conjecture that
$H^{1,*}_{{\A}_{m}}(P_{v})$ can be described as
follows\footnote{Conjectures \ref{HAm,1 for a polygon} and
\ref{homology polygon Am} have now been proved in \cite{P05}.}.
\begin{conjecture} $H_{\A_m}^1$ for a polygon\
\begin{enumerate}
\item[(Odd)] $Tor(H^{1,*}_{{\A}_{m}}(P_{2g+1}))=
H^{1,gm}_{{\A}_{m}}(P_{2g+1})=\Z_m$.\\
The Poincar\'{e} polynomial of $H^{1,*}_{{\A}_{m}}(P_{2g+1})$ is
equal to\\
 $tq^{(g-1)m}(q+q^2+...+q^{m-1})$.
\item[(Even)] $H^{1,*}_{{\A}_{m}}(P_{2g+2})$ is a free abelian
group with the Poincar\'{e} polynomial equal to
$tq^{gm}(q+q^2+...+q^{m-1})$.
\end{enumerate}
\label{HAm,1 for a polygon}
\end{conjecture}

We checked the conjecture in  a few other specific cases in addition
to $m=2$. In particular, we verified that $H_{\A_3}^{1,6}(P_5)=\Z_3$
and that $H_{\A_3}^{1}(P_4)=\Z\{4\} \oplus \Z\{5\}$.

\bigskip

This conjecture can be used to describe in full the cohomology of
$P_v$ via Corollary \ref{Polygon in terms of H1}. Note that for
$m=2$ we get Theorem \ref{Homology polygon A2} from Section
\ref{section algebra is A2}.

\begin{conjecture}All cohomology groups $H_{\A_m}^*$ for a polygon:\
\label{homology polygon Am}
\begin{enumerate}
\item[(Odd)] For $v=2g+1$ we have: \\
$Tor(H^{*,*}_{{\A}_{m}}(P_{2g+1}))=
 H^{v-2,m}_{{\A}_{m}}(P_{2g+1})\oplus
H^{v-4,2m}_{{\A}_{m}}(P_{2g+1})\oplus ... \oplus
H^{1,gm}_{{\A}_{m}}(P_{2g+1})$ with each summand
isomorphic to $\Z_m$.\\
The Poincar\'{e} polynomial of $H^{*,*}_{{\A}_{m}}(P_{2g+1})$ is
equal to\\
$ (q+...+q^{m-1})^v +$ \\
$(q+...+q^{m-1})(t^{v-2}+ (t^{v-3}+t^{v-4}q^m + (t^2+t)q^{m(g-1)})$.
\item[(Even)]
For $v=2g+2$ we have: \\
$Tor(H^{*,*}_{{\A}_{m}}(P_{2g+2}))=
H^{v-2,m}_{{\A}_{m}}(P_{2g+2})\oplus
H^{v-4,2m}_{{\A}_{m}}(P_{2g+2})\oplus ... \oplus
H^{2,gm}_{{\A}_{m}}(P_{2g+2})$ with each summand
isomorphic to $\Z_m$.\\
The Poincar\'{e} polynomial of $H^{*,*}_{{\A}_{m}}(P_{2g+2})$ is
equal to\\
$(q+...+q^{m-1})^v + q^{m(v/2)-1)}(q+...+q^{m-1})+$ \\
$(q+...+q^{m-1})(t^{v-2}+ (t^{v-3}+t^{v-4})q^m +
(t^3+t^2)q^{m(g-1)})+ tq^{mg})$.
\end{enumerate}
\end{conjecture}

 The knowledge of torsion in cohomology of the triangle
 allows us, in some cases, to deduce
existence of torsion in other graphs.

\bigskip

A $v$-gon with diagonals is a graph obtain from a polygon by adding
edges between some vertices of the polygon such that it is possible
to draw all the diagonals inside the polygon without them
intersecting.

\begin{proposition}\label{$v$-gon with diagonals}
If $G$ is a $v$-gon with diagonals and $\A =\Z1 \oplus \A'$ is an
algebra satisfying Assumption \ref{assumptions on algebras}, then
$H^{v-2,*}_{\A}(G)= H^{1,*}_{\A}(P_3)$. \\
In particular, $H^{v-2,m}_{\A_m}(G)= \Z_m$.
\end{proposition}
\noindent \emph{Proof.}\ \ The proposition follows by induction on
the number of edges from the following two lemmas which allow the
reduction of the cohomology of $G$ to the one of the triangle.

\begin{lemma}
\label{lemma both endpoints degree two} Let $e$ be an edge of $G$
such that both endpoints of $e$ have degree two and let $v=v_G
>3$. Then
$$ H_{\A}^{v_{G/e}-2,*}(G/e) = H_{\A}^{v_G-2,*}(G).$$
\end{lemma}
\begin{proof}
Since $G-e$ has two pendant edges, by Corollary \ref{pendant edge
vanishing} part (i), $H_{\A}^{v-3,*}(G-e)=0$ and
$H_{\A}^{v-2,*}(G-e)=0$.


It follows that the long exact sequence of cohomology for $G$ at the
edge $e$
 reduces in our case to
$$ 0=H_{\A}^{v-3,*}(G-e) \to H_{\A}^{v-3,*}(G/e) \to H_{\A}^{v-2,*}(G) \to
H_{\A}^{v-2,*}(G-e)=0$$ and the lemma follows.
\end{proof}
\begin{lemma}
\label{e on a triangle and vertex opposite  has degree two} Let
$e$ be an edge of $G$ such that $e$ lies on a triangle in $G$ and
the vertex on the triangle opposite to $e$ has degree two. If
$v=v(G)
> 3$ then
$$ H_{\A}^{v-2,*}(G-e) = H_{\A}^{v-2,*}(G).$$
\end{lemma}
\begin{proof}
By Remark \ref{Remark6}, $H_{\A}^{v_{G/e}-1,*}(G/e)\cong 0$. Since
$v_{G/e}=v_G-1$, $H_{\A}^{v_{G}-2,*}(G/e)\cong 0$. Also, after
replacing the double edge of $G/e$ by a single edge (we know that
this doesn't change the cohomology groups \cite{HR05}), $G/e$ has
one pendant edge, so by Corollary \ref{pendant edge vanishing} part
(i), $H_{\A}^{v_{G/e}-2,*}(G/e)\cong 0$ i.e.
$H_{\A}^{v_G-3,*}(G/e)\cong 0$. The long exact sequence of
cohomology reduces in this case to
$$ 0=H_{\A}^{v_G-3,*}(G/e) \to H_{\A}^{v_G-2,*}(G) \to H_{\A}^{v_G-2,*}(G-e) \to
H_{\A}^{v_G-2,*}(G/e)=0$$ and the lemma follows.
\end{proof}

If we relax the assumptions of Lemma \ref{lemma both endpoints
degree two} we get the following weaker result.

\begin{proposition}
\label{lemma one endpoint degree two} Let $e$ be an edge of $G$ such
that one endpoint of $e$ has degree two and let $v=v_G \geq 3$.\\
Then there in an epimorphism $\alpha: H_{\A}^{v_{G/e}-2,*}(G/e)
\twoheadrightarrow H_{\A}^{v_G-2,*}(G).$
\end{proposition}
The proof of the proposition is that same as that of Lemma
\ref{lemma both endpoints degree two}, except that we use only the
right side of the exact sequence from the proof. Notice that we
cannot expect $\alpha$ to be always an isomorphism; for example for
$G=P_3$ the group $H_{\A}^{0,*}(G/e)$ is free while
$H_{\A}^{1,*}(G)$ has a torsion. Proposition \ref{lemma one endpoint
degree two} is not even sufficient to guarantee that if
$H_{\A}^{v_{G/e}-2,*}(G/e)$ has a torsion then $H_{\A}^{v_G-2,*}(G)$
has a torsion. It can be used however to show that for a large class
of graphs (corresponding in Tait translation to nontrivial, prime
algebraic link diagrams) the group $H_{\A}^{v_G-2,*}(G)$ is
nontrivial and dominates the first homology group of the triangle.

\begin{definition}
A {\it $P_n$-series-parallel graph} is a graph which can be obtained
from the polygon $P_n$ by the finite number of the following
operations: replacing an edge by two edges or subdividing an edge.
\end{definition}

\begin{proposition}
\label{Pn series-parallel}
 If a graph $G$ is
$P_n$-series-parallel then there is an epimorphism
$H^{n-2}_{\A}(P_n)\twoheadrightarrow H^{v-2}_{\A}(G)$, or
equivalently $H^1_{\A}(P_3)\twoheadrightarrow H^{v-2}_{\A}(G)$ since
the two groups on the left are isomorphic by Corollary \ref{Polygon
in terms of H1}.

\end{proposition}
\begin{proof}
We note that adding a parallel edge does not change the cohomology
groups \cite{HR05}, and subdividing an edge induces an epimorphism
(Proposition \ref{lemma one endpoint degree two}). Therefore the
result follows by an induction on the number of edges.
\end{proof}

Since any  $P_{n}$-series-parallel graph is also a
$P_{3}$-series-parallel graph, this combined with Proposition
\ref{thm Am P3} implies that if a graph $G$ is
$P_{n}$-series-parallel for some n, $H^{v-2,m}_{\A_m}(G) $ is equal
to $\Z_{m'}$ or $0$, where $m'$ is a divisor of $m$.

 We can use the
following proposition to show that some class of graphs, including
the complete graph $K_4$, has torsion in $H_{\A_m}^{v-2,m}(G)$.

\begin{proposition}
\label{Hv-2 G torsion when G/e finite and G-e torsion} Let $e$ be
an edge of a graph $G$. If $H_{\A}^{v-3,j}(G/e)$ is finite,
possibly trivial, and $H_{\A}^{v-2,j}(G-e)$ has a torsion, then
$H_{\A}^{v-2,j}(G)$ has a torsion.
\end{proposition}
\begin{proof}
It follows from the long exact sequence of graph cohomology:
$$ \to H_{\A}^{v-3,j}(G/e) \to H_{\A}^{v-2,j}(G) \to H_{\A}^{v-2,j}(G-e)
\to H_{\A}^{v-2,j}(G/e)=0$$
\end{proof}

To illustrate Proposition \ref{Hv-2 G torsion when G/e finite and
G-e torsion}, consider the complete graph $K_4$ and any edge, $e$ of
it. We know by Proposition \ref{$v$-gon with diagonals} that
$H^{2,m}_{\A_m}(K_4-e)=\Z_m$ and that $H^{1,m}_{\A_m}(K_4/e)= \Z_m$
therefore we have an exact sequence $ \cdots \to \Z_m \to
H^{2,m}_{\A_m}(K_4) \to \Z_m \to 0$. We conclude from this that
$H^{2,m}_{\A_m}(K_4)$ is a torsion group (extension of $\Z_{m'}$ by
$\Z_m$, where $m'$ is a divisor of $m$) and contains element of
order $m$. Recall that a group $H$ is an extension of $\Z_{m'}$ by
$\Z_m$ means that it satisfies an exact sequence of the form $0
\rightarrow \Z_{m'} \rightarrow H \rightarrow  \Z_{m}\rightarrow 0$.

\begin{figure}[h]
\begin{center}
\scalebox{.4}{\includegraphics{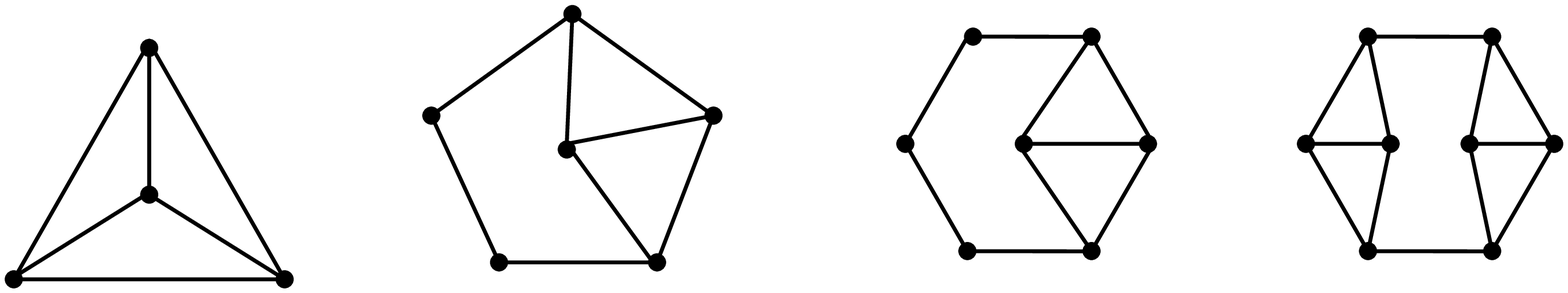}}
\caption{\textit{Graphs with torsion in
$H^{v_G-2,m}_{\A_m}(G)$.}\label{finHomgraphs}}
\end{center}
\end{figure}

The graph $K_4$ is one of the series of graphs for which we can
prove, in a similar manner, that $H^{v_G-2,m}_{\A_m}(G)$ has
nontrivial torsion. Some of these graphs are illustrated in Figure
\ref{finHomgraphs}.

Notice, that in our examples $H_{\A_m}^{v_G - 2}(G)$ had nontrivial
torsion. It is not always the case, however, at least for graphs
which are not 2-connected. Consider, for example, the one vertex
product $P_3*P_3=$
\parbox{1.1cm}{\psfig{figure=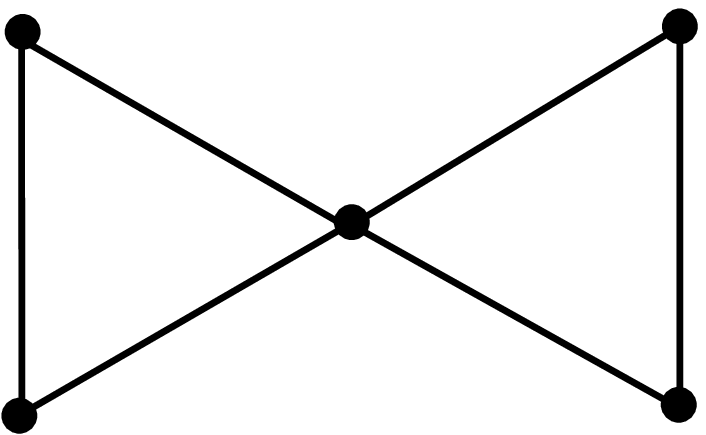,height=0.5cm}}.
Then $H_{\A}^{3}(P_3*P_3) = 0$. This is the consequence of the
following three term part of the
 deleting-contracting long exact sequence:
$$ 0=
H_{\A}^{2}(\parbox{0.5cm}{\psfig{figure=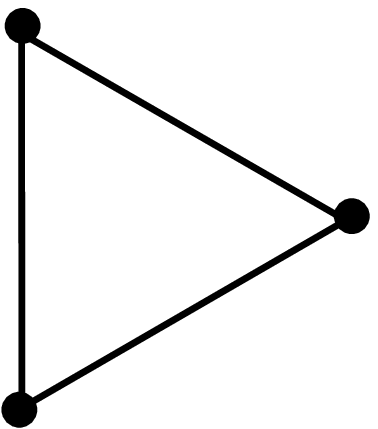,height=0.5cm}})
\otimes \A' =H_{\A}^{2}(
\parbox{1.0cm}{\psfig{figure=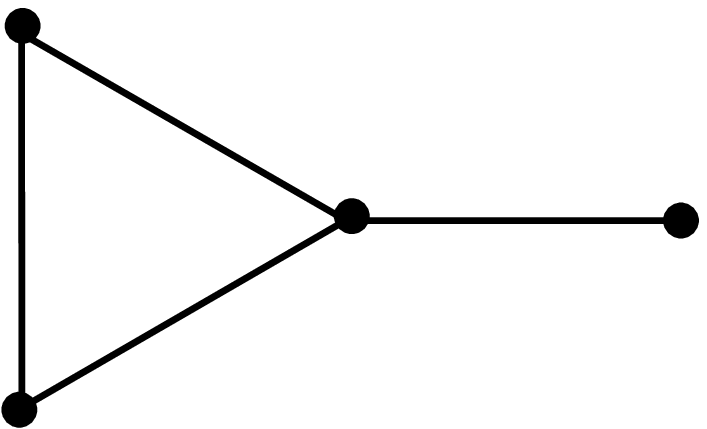,height=0.5cm}})=
H_{\A}^{2}(
\parbox{1.0cm}{\psfig{figure=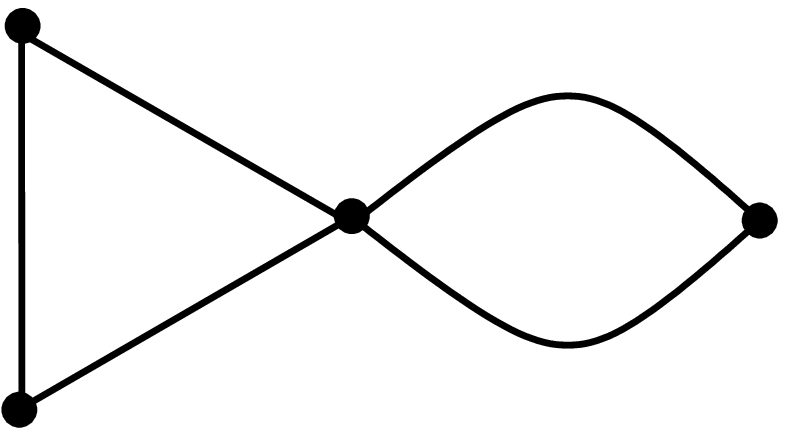,height=0.5cm}})
\overset{\alpha_*}{\rightarrow}
H_{\A}^{3}(\parbox{1.0cm}{\psfig{figure=Gbutterfly.eps,height=0.5cm}})
\overset{\beta_*}{\rightarrow} $$
$$H_{\A}^{3}(
\parbox{1.0cm}{\psfig{figure=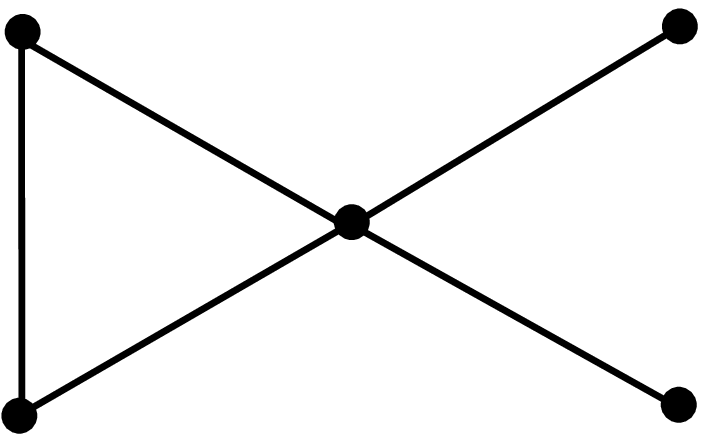,height=0.5cm}})=
H_{\A}^{3}(\parbox{0.5cm}{\psfig{figure=Gtriangle.eps,height=0.5cm}})
\otimes \A' \otimes \A' = 0.$$

In a broader context one can try to generalize results obtained for
$m=2$ from Section \ref{section algebra is A2}. One would conjecture
that if a loopless graph $G$ contains a triangle than
$H^{1,*}_{\A_m}(G)$ contains a $\Z_m$-torsion and if a loopless
graph $G$ contains a square than $H^{2,*}_{\A_m}(G)$ contains a
$\Z_m$-torsion.

\medskip

\begin{figure}[h]
\begin{center}
\scalebox{.3}{\includegraphics{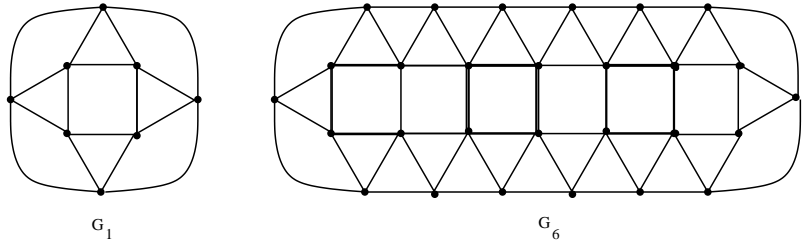}}
\caption{\textit{Family of graphs $G_i$}, where $i$ is the number of
squares.} \label{FamilyGi}
\end{center}
\end{figure}

A. Shumakovitch has conjectured that for alternating links the
torsion in Khovanov homology can have only elements of order $2$. We
conjecture analogously that for the algebra $\A_2$ the torsion part
of $H^{*,*}_{\A_2}(G)$ can  have only elements of order $2$.
However, as computed in \cite{PPS06}, $H^{1,5}_{A_3}(K_4) = \Z_3^2
\oplus \Z_6 \oplus \Z^2$ so it is not true that torsion elements of
$H^{**}_{\A_3}(G)$ have order a power of three. Futhermore, in
\cite{PPS06}, the family of plane graphs $G_1$, $G_2$,...,$G_k$ (see
Figure \ref{FamilyGi}) is analysed and there is a strong indication
that for any $n$ there is a graph with torsion $\Z_n$. In
particular, $H^{1,13}_{\A_3}(G_1)= \Z_4 \oplus \Z_3^7 \oplus
\Z^{12}$, $H^{1,21}_{\A_3}(G_2)= \Z_{18}\oplus \Z_3^{10}\oplus
\Z^{15}$, $H^{1,29}_{\A_3}(G_3)= \Z_{8}\oplus \Z_3^{15}\oplus
\Z^{20}$, $H^{1,37}_{\A_3}(G_4)= \Z_{10}\oplus \Z_3^{19}\oplus
\Z^{25}$, $H^{1,45}_{\A_3}(G_5)= \Z_{4}\oplus \Z_9 \oplus
\Z_3^{22}\oplus \Z^{30}$, $H^{1,53}_{\A_3}(G_6)= \Z_{14}\oplus
\Z_3^{27}\oplus \Z^{35}$, and $H^{1,61}_{\A_3}(G_7)= \Z_{16}\oplus
\Z_3^{31}\oplus \Z^{40}$.


\subsection{Deformation of the algebra $\A_m$} \label{Deformation5.3}

Khovanov and Rozansky constructed $sl(m)$ homology for links
\cite{KR04}. The underlying algebra in this case was
$\A_m=\Z[x]/(x^m)$. Then Gornik considered deformations of this
$sl(m)$ homology using as the underlying algebra
$\A_{p(x)}=\Z[x]/p(x)$, where $p$ is a polynomial.

Here, we prove that\footnote{The algebra $\A_{p(x)}$ is not
necessary free, for example $\A_{2x}=\Z \oplus \Z_2 \oplus
\Z_2\oplus ...$ and $\A_{2x+1}$ is a torsion free but not free
group. In fact, $\A_{2x+1} =\Z[\frac{1}{2}]$. To be sure that
$\A_{p(x)}$ is free we can assume that $p(x)$ is monic, but, as
mentioned in the introduction, many basic properties of graph
cohomology holds without the assumption,
 in fact no restrictions on $p(x)$ are needed in a proof
of Theorem \ref{H1 P3 Ap(x)}.}:
\begin{theorem}\
\label{H1 P3 Ap(x)}
\begin{enumerate}
\item[(i)] $H^{1}_{\A_{p(x)}}(P_3)=\Z[x]/(p(x),p'(x))$, where
$p'(x)$ is the derivative of $p(x)$ and $(p(x),p'(x))$ is the ideal
generated by $p(x)$ and $p'(x)$. In particular, \emph{rank}$\:(
H^{1}_{\A_{p(x)}}(P_3))=\deg_\mathbb{Q}(gcd(p(x),p'(x))$.
Equivalently, \emph{rank}$\:( H^{1}_{\A_{p(x)}}(P_3))=\sum_{x_i}
(mul(x_i)-1)$, where the sum is taken over all roots of $p(x)$ and
$mul(x_i)$ is the multiplicity of of the root $x_i$.
\item[(ii)]
 $H^{0}_{\A_{p(x)}}(P_3)=\Z^{m(m-1)(m-2)+\deg_{\mathbb{Q}}(gcd(p(x),p'(x))}$,
where $m=\deg(p(x))$ is the degree of $p(x)$.
\item[(iii)] More generally, consider the algebra $\A=\Z[x]/\mathbb
I$ where $\mathbb I$ is an ideal in $\Z[x]$. Let $\mathbb I'$ be the
ideal generated by derivatives of elements in $\mathbb I$. Note that
$\mathbb I \subset \mathbb I'$ because for any $p(x)\in \mathbb I$
we have $p(x)= (xp(x))' -xp'(x)$. Then $H^{1}_{\Z[x]/\mathbb I}=
\Z[x]/\mathbb I'$.
\end{enumerate}
\end{theorem}

Note that for $p(x)=x^m$, Theorem \ref{H1 P3 Ap(x)} agrees with
Theorem \ref{thm Am P3} because $\Z[x]/(x^m,mx^{m-1})$, as a group,
is isomorphic to $\Z^{m-1} \oplus \Z_m$.

For $p(x)= x^2-bx-a$ we recover a result from \cite{HR05} because
$\Z[x]/(x^2-bx-a, 2x-b) = \Z[x]/(x^2-bx-a, 2x-b, bx+2a)$ which as a
group is isomorphic to $\{1,x\ | \ 2x-b, bx+2a \}$ and finally to
$\Z \oplus \Z_2$ for $b^2+4a=0$, $\Z_{|b^2+4a|}$ for $b^2+4a\neq 0$
and $b$ odd and $\Z_2 \oplus \Z_{|b^2+4a|/2}$ for $b^2+4a\neq 0$ and
$b$ even.

For $p(x)=x^m-1$ we get from Theorem \ref{H1 P3 Ap(x)} that
$H^{1}_{\A_{p(x)}}(P_3) = \Z[x]/(x^m-1, mx^m)$ which is isomorphic
to $\Z_m^m$. In fact, for $p(x)=x^m-1$ we conjecture more generally
that $H_{\A_{x^m -1}}^1(P_v)$ is equal to $0$ for even $v$ ($v\geq
2$) and it is equal to $\Z_m^m$ for odd $v$.

\begin{proof}[Proof of Theorem \ref{H1 P3 Ap(x)}]

Theorem \ref{H1 P3 Ap(x)} follows from Lemma \ref{main lemma} by
carefully analysing relations caused by dividing $\Z[x]$ by $p(x)$
or $\mathbb I$. To simplify the proof we introduce several useful
general notions enriching the structure of $H^*_{\A}(G)$ for any
$\A$ and $G$.

Recall that an enhanced state $S$ of $C^*_{\A}(G)$ is a spanning
subgraph of $G$ with an element of $\A$ associated to each connected
component. A typical element of $C^*_{\A}(G)$ is a linear
combination of enhanced states with coefficients in $\Z$.

\begin{proposition}
\label{Homology groups are A-modules} For a graph $G$ and a chosen
base vertex, $v_1$, the group $H^*_{\A}(G)$
is an $\A$-module as follows: \\
$\A$ acts on $C^*_{\A}(G)$ the following way: For any $a \in \A$ and
any enhanced state $S \in C^*_{\A}(G)$, $a\cdot S$ is obtained from
$S$ by multiplying by $a$ the weight of the component
 of $[G:s]$ containing $v_1$.  The action commutes with
boundary map $d$, therefore $H^*_{\A}(G)$ is an $\A$-module.
\end{proposition}
\begin{proof}
The proof is left to the reader.
\end{proof}

We will also need the following useful tool.

\begin{proposition}\label{Let I be an ideal}
Let ${\mathbb I}$ be an ideal in the algebra $\A$. For a given
graph $G$ consider the short exact sequence of chain complexes (of
$\A$-modules):
$$           0 \to  \ker(P) \overset{I}{\rightarrow}  C^*_{\A}(G)
        \overset{P}{\rightarrow} C^*_{\A/\mathbb I}(G) \to 0,$$
where $P$ is induced by the projection $\A \to \A/\mathbb I$ and
$I$ is the inclusion map. This short exact sequence yields a long
exact sequence of cohomology groups (and $\A$-modules):
$$ \ldots \to H^i_{\A,\mathbb I}(G) \overset{I_*}{\rightarrow} H^i_{\A}(G)
\overset{P_*}{\rightarrow} H^i_{\A/{\mathbb I}}(G)
\overset{\partial}{\rightarrow} H^{i+1}_{\A,\mathbb I}(G)\to
\ldots$$ We call the group ($\A$-module) $H^i_{\A,\mathbb I}(G)$ the
\emph{relative graph cohomology of $G$}.
\end{proposition}


We use now the exact sequence of Proposition \ref{Let I be an
ideal}
 to prove Theorem \ref{H1 P3 Ap(x)}. For a triangle
$P_3$ and ${\mathbb I} \subset \Z[x]$. We have
$$ \to H^1_{\Z[x],\mathbb I}(P_3) \overset{I_*}{\rightarrow} H^1_{\Z[x]}(P_3)
\overset{P_*}{\rightarrow} H^1_{\Z[x]/{\mathbb I}}(P_3)
\overset{\partial}{\rightarrow} H^{2}_{\Z[x],\mathbb I}(P_3)$$

In order to use this exact sequence to prove Theorem \ref{H1 P3
Ap(x)} we have to understand $H^2_{\Z[x],\mathbb I}(P_3)$ and
$H^1_{\Z[x],\mathbb I}(P_3)$. The ring $\Z[x]$ is a N{\oe}terian
ring so $\mathbb I$ is finitely generated, say by polynomials
$p_1(x),...,p_k(x)$. We denote by $(p_j(x))$ the principal ideal
generated by $p_j(x)$.

\begin{lemma} \ \\
\begin{enumerate}
\item[(i)] $H^2_{\Z[x],\mathbb I}(P_3) = 0$, \item[(ii)]
$H^1_{\Z[x],\mathbb I}(P_3)$ is generated by the classes of
 elements of of the form (notation is explained in
Figure \ref{C1C2C3basis_tr}):\\
$f_1^{p_j(x)x^{\ell},1} + f_2^{1,p_j(x)x^{\ell}} +
f_3^{p_j(x)x^{\ell},1}$ for any $j$ and $\ell$, and
$f_1^{u,w}-f_1^{p_j(x)x^{\ell},1}$ with $uw=p_j(x)x^{\ell}$ and $u$
or $w$ in $(p_j(x))$, and their rotations
$\tau(f_1^{u,w}-f_1^{p_j(x)x^{\ell},1})=
f_2^{w,u}-f_2^{1,p_j(x)x^{\ell}}$ and
$\tau^2(f_1^{u,w}-f_1^{p_j(x)x^{\ell},1})=
f_3^{u,w}-f_3^{p_j(x)x^{\ell},1}$.
\end{enumerate}
\label{H2-H1}
\end{lemma}
\begin{proof}
(i)In the relative cochain complex the group
$C^2_{\Z[x],\mathbb{I}}$ is generated by elements $f^{u}_{1,2},
f^{u}_{2,3}, f^{u}_{1,3}$
 with $u\in \mathbb{I}$; see Figure
\ref{C1C2C3basis_tr}. $\ker(d^2_{\Z[x],\mathbb{I}})=
\ker(C^2_{\Z[x],\mathbb{I}} \to C^3_{\Z[x],\mathbb{I}})$ is
generated by elements of the form $f^{u}_{1,2} - f^{u}_{2,3}$ and
$f^{u}_{1,3} + f^{u}_{2,3}$, $u\in \mathbb{I}$. These elements are
in the image of $C^1_{\Z[x],\mathbb{I}}$; namely, they are equal
to $d^1(f_2^{1,u})$ and $d^1(f_3^{u,1})$ respectively; $u\subset
\mathbb{I}$ in our considerations. Therefore $H^{2}_{\Z[x],\mathbb
I}(P_3)=0$.

(ii) In fact, $d^1_{\Z[x],\mathbb{I}}:
\mbox{Span}(f_2^{1,u},f_3^{u,1}) \to \ker(d^2_{\Z[x],\mathbb{I}})$
is an isomorphism. The $\Z$-module Span$(f_2^{1,u},f_3^{u,1})$ can
be written equivalently as
Span$(f_2^{1,p_j(x)x^{\ell}},f_3^{p_j(x)x^{\ell},1})$ for any $j$
and $\ell$. The cochain group $C^1_{\Z[x],\mathbb I}(P_3)$ has a
natural generating set $f_i^{x^{\ell_1},p_j(x)x^{\ell_2}}$, and
$f_i^{p_j(x)x^{\ell_1},x^{\ell_2}}$ for $i=1,2,3$ and any $j$,
$\ell_1$ and $\ell_2$. We can decompose the group as
$C^1_{\Z[x],\mathbb I}(P_3) = \ker(d^2_{\Z[x],\mathbb{I}}) \oplus
\mbox{Span}(f_2^{1,p_j(x)x^{\ell}},f_3^{p_j(x)x^{\ell},1})$.
Therefore $\ker(d^2_{\Z[x],\mathbb{I}})$ is generated by
$f_1^{p_j(x)x^{\ell},1} + f_2^{1,p_j(x)x^{\ell}} +
f_3^{p_j(x)x^{\ell},1}$ for any $j$ and $\ell$, and
$f_1^{u,w}-f_1^{p_j(x)x^{\ell},1}$ with $uw=p_j(x)x^{\ell}$ and $u$
or $w$ in $(p_j(x))$, and $f_2^{w,u}-f_2^{1,p_j(x)x^{\ell}}$ and
$f_3^{u,w}-f_3^{p_j(x)x^{\ell},1}$. Classes of generators of
$\ker(d^2_{\Z[x],\mathbb{I}})$ generate $H^1_{\Z[x],\mathbb
I}(P_3)$.
\end{proof}


\begin{figure}[h]
\begin{center}
\scalebox{.45}{\includegraphics{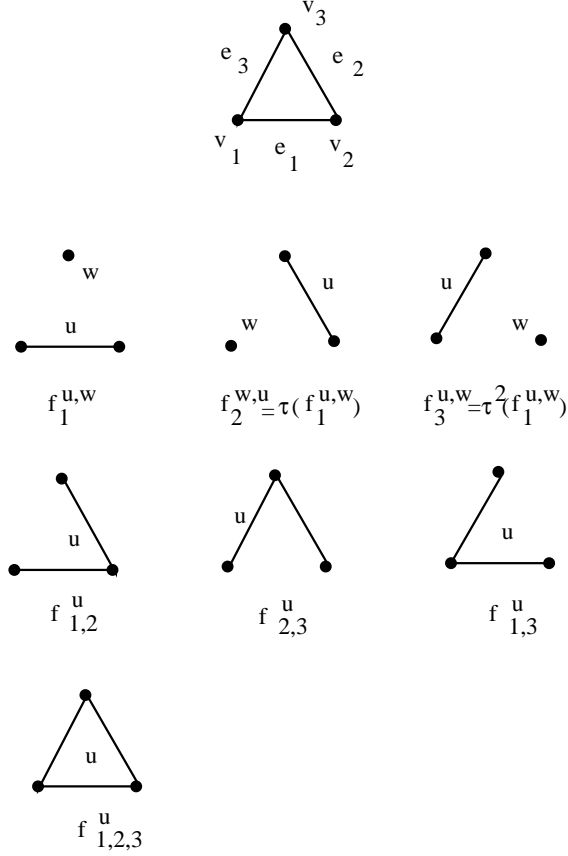}}
\caption{\textit{Elements of the cochain group $C^{i}_{\A/{\mathbb
I}}(P_3)$;\ $i>0$, $u,v \in \A$. In $f^{u,v}_i$, $u$ indicates the
color of the component containing $v_1$.} \label{C1C2C3basis_tr} }
\end{center}
\end{figure}

We can complete now the proof of Theorem \ref{H1 P3 Ap(x)}. It
follows from the fact that $P_*$ is an epimorphism that
$H^1_{\Z[x]/\mathbb{I}}(P_3)$ is generated by the image of the basis
of $H^1_{\Z [x]}(P_3)$ that is  $e=e_1 \equiv f_1^{1,x}-f_1^{x,1}= $
\parbox{3.3cm}{\psfig{figure=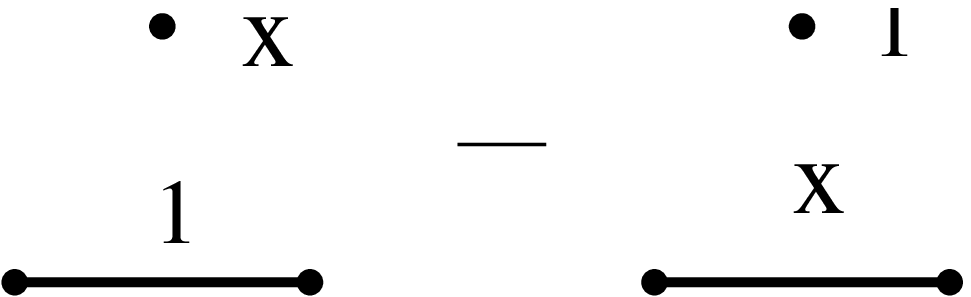,height=0.8cm}},\\ \ \\
$e_2=xe \equiv f_1^{x,x}-f_1^{x^2,1}$, ..., $e_j \equiv
x^{j-1}e=f_1^{x^{j-1},x}-f_1^{x^j,1}$,... (notation
of Figure \ref{C1C2C3basis_tr}).\\
In order to find presentation of  $H^1_{\Z[x]/{\mathbb I}}(P_3)$ we
have to analyze relations given by Im$(I_*)$. By Lemma \ref{H2-H1}
(ii) we know the generating set of $H^1_{\Z[x],\mathbb I}(P_3)$, so
we will compute their images under $I_*$. First we check that
$I_*(f_1^{p_j(x),x} - f_1^{xp_j(x),1})=
p_j(x)(f_1^{1,x}-f_1^{1,x})=p_j(x)e$. Then, by  Lemma \ref{main
lemma}, assuming that $p_j(x)=\sum_{i=0}^m a_ix^i$, we have
$I_*(f_1^{1,p_j(x)}- f_1^{xp_j(x),1})= f_1^{1,\sum a_ix^i} -
f_1^{\sum a_ix^i,1} = \sum_{i=0}^m a_i (f_1^{1,x^i} - f_1^{x^i,1})
\equiv $ $\sum_{i=0}^m ia_i (f_1^{x^{i-1},x} - f_1^{x^i,1})=
\sum_{i=0}^m ia_ie_i =
 \sum_{i=0}^m ia_ix^{i-1}e = p_j'(x)e$.\\
Polynomials $p_j(x)$ and $p_j'(x)$ (any $j$) generate the ideal
$\mathbb I'$. Therefore $I_*(H^1_{\Z[x],\mathbb I}(P_3))$ contains
$I'e$. We will argue now that $I_*(H^1_{\Z[x],\mathbb I}(P_3))=
\mathbb I'e$; therefore $H^{1}_{\Z[x]/\mathbb I}= \Z[x]/\mathbb
I'$. Now we perform calculations checking images of generators of
$H^1_{\Z[x],\mathbb I}(P_3)$ under $I_*$:\\
$I_*(f_1^{x^{k_1}p_j(x),x^{k_2}} - f_1^{x^{k_1}x^{k_2}p_j(x),1}) =
x^{k_1}p_j(x)(f_1^{1,x^{k_2}} - f_1^{x^{k_2},1}) =
kx^{k_1+k_2-1}p_j(x)e
\in \mathbb I'e$. \\
$I_*(f_1^{x^{k_1},x^{k_2}p_j(x)} - f_1^{x^{k_1}x^{k_2}p_j(x),1}) =
x^{k_1}(f_1^{1,x^{k_2}p_j(x)} - f_1^{x^{k_2}p_j(x),1}) =
x^{k_1}\sum_{i=0}^m (i+k_2)x^{i+k_2-1}e = kx^{k_1}x^{k_2-1}p_j(x)e
+
x^{k_1+k_2}p'_j(x)e \in \mathbb I'e$. \\
Similarly we deal with with $\tau$ rotations of this elements. We
are left to check $I_*(f_1^{x^kp_j(x),1} + f_2^{1,x^kp_j(x)}
+f_3^{x^kp_j(x),1}) = x^kp_j(x)(f_1^{x,1} + f_2^{x,1} + f_3^{x,1})
-(f_2^{x^kp_j(x),1} - f_2^{1,x^kp_j(x)})= x^kp_j(x)(f_1^{x,1} +
f_2^{x,1} + f_3^{x,1}) - kx^{k-1}p_j(x)e-x^kp_j'(x)e \in \mathbb
I'e$. The proof of Theorem \ref{H1 P3 Ap(x)} (iii) and (i) is
completed. Part (ii) follows by considering the chromatic polynomial
of the triangle $\lambda(\lambda-1)(\lambda -2)$. This completes the
proof of Theorem \ref{H1 P3 Ap(x)}.
\end{proof}

\bigskip

Potentially the results and ideas used in Theorem \ref{Let $D$ be
the diagram} can be used to compute partially the Khovanov-Rozansky
$sl(n)$ homology and its deformations as defined by Gornik. In
particular, Theorem \ref{thm Am P3} potentially predicts
Khovanov-Rozansky $sl(n)$ homology for the trefoil. Conjecture
\ref{homology polygon Am} potentially predicts Khovanov-Rozansky
$sl(n)$ homology for the $(2,n)$ torus links. In particular, we can
predict some torsion in Khovanov-Rozansky $sl(n)$ homology in these
cases.

Our calculations also confirm the connection between graph
cohomology and the Hochschild homology of the underlying algebra.
Furthermore, we also speculate a connection between our graph
cohomology of a symmetric graph and the Cones cyclic homology of
algebra \cite{P05}.

\begin{small}
\noindent Department of Mathematics, \\
 The George Washington
University,\\
1922 F street NW, \\
Washington, DC 20052, U. S. A. \\
Email: lhelmeg@yahoo.com, przytyck@gwu.edu, rong@gwu.edu.
\end{small}

\end{document}